\documentclass[final,3p,times]{elsarticle}
\usepackage[fleqn]{amsmath}
\usepackage{float}
\usepackage{amssymb}
\usepackage{amsthm,mathrsfs}
\usepackage{lineno}
\usepackage{amsfonts}
\usepackage{epsfig}
\usepackage{color}
\usepackage{setspace}
\usepackage{latexsym}
\biboptions{numbers,sort&compress}
\usepackage[colorlinks,citecolor=blue,urlcolor=blue,linkcolor=blue]{hyperref}
\usepackage{subfigure,graphicx,epstopdf}
\usepackage{exscale,relsize}
\allowdisplaybreaks

\begin{document}
\begin{frontmatter}

\title{Existence of one class of global strong solution to the Cauchy problem for the three-dimensional incompressible Boussinesq system}

\author[labela]{Rulv Li}\ead{rulvl@emails.bjut.edu.cn}
\author[labela]{Shu Wang\corref{cor1}}\ead{wangshu@bjut.edu,cn}

\address[labela]{School of Mathematics, Statistic and Mechanics, Beijing University Of Technology,
Beijing, 100124, China.} \cortext[cor1]{Corresponding author}

\begin{abstract}
In this paper, we prove the the global existence of strong solutions to the three dimensional incompressible Boussinesq system with some special solenoidal initial data. In particular, these solutions can be expressed into the Fourier series.
\end{abstract}


\end{frontmatter}


\section{Introduction}
In this paper we study the existence of one class of global strong solution to the following Cauchy problem for three-dimensional incompressible Boussinesq system
\begin{eqnarray}\label{eq1}
\begin{cases}
\partial_tu+u\cdot\nabla u-\mu\Delta u+\nabla \pi=\rho e_3,~~~&(t,x)\in(0,T]\times\mathbb{R}^3,\\
\partial_t\rho+u\cdot\nabla\rho-\nu\Delta\rho=0,~~~&(t,x)\in(0,T]\times\mathbb{R}^3,\\
div~u=0,~~~&(t,x)\in(0,T]\times\mathbb{R}^3,\\
u(0,x)=u_0(x),\rho(0,x)=\rho_0(x)~~~&x\in\mathbb{R}^3.
\end{cases}
\end{eqnarray}
Here the unknown functions $u,\pi,\rho$ are the velocity vector field, the scalar pressure and the temperature respectively, the constant $T>0$ is any given time, $e_3=(0,0,1)$ id given unit vector and the initial data $u_0(x), \rho_0(x)$ are given smooth functions with $div~u_0(x)=0$, positive number $\mu,\nu$ is viscosity and thermal diffusivity.

Thanks to the physical and mathematical significances, the global well-posedness of \eqref{eq1} becomes the main focus in fluid dynamics. In the case of two-dimensions, the unique global solution for $\mu>0$ and $\nu>0$ can be obtained by the classical energy methods (see \cite{6}). For the partial diffusion case, Chae \cite{8} or Hou-Li \cite{15} proved the global regularity for the Boussinesq system in $H^s (s\geq2)$ when $\mu>0$, $\nu=0$ as well as $\mu=0$, $\nu>0$. Later, Cao-Wu \cite{7} proved the global existence of classical solutions to the anisotropic Boussinesq equations and first found that the growth of the $L^p$-norm of the vertical velocity is not faster than
$\sqrt{r\ln r}$ at infinity. We refer to \cite{2,4,5} for more relevant results. For the fractional diffusion case, Hmidi-Keraani-Rousset \cite{14} made full use of the cancellation condition and showed the existence and uniqueness of the solution in critical framework. Miao-Xue \cite{20} proved the global well-posedness of the Boussinesq equations with fractional viscosity and thermal diffusion under the assumption that the
fractional powers satisfy a mild condition. Xu-Ye \cite{25,26} established the global regularity of the smooth solutions for a new range of fractional powers of the diffusion.

For the three-dimensional Boussinesq equations, Danchin-Paicu \cite{10} showed the global existence of weak solution for $L^2$-data and the global well-posedness for small initial data, and in \cite{9} also obtained an existence and uniqueness result for small initial data belonging to some critical Lorentz spaces. As for three-dimensional incompressible Navier-Stokes equations, the question on global existence of smooth solutions to three-dimensional incompressible Boussinesq equations is still open. In the case of no swirl, inspired by the global results for Euler equations such as \cite{24}, Abidi-Hmidi-Sahbi \cite{1} showed the global well-posedness with $supp ~\rho_0$ does not intersect the axis $oz$ for $\mu>0,\nu=0$ and Hmidi-Rousset \cite{13} also obtained the same result for $\mu=0,\nu>0$. Miao-Zheng \cite{21,22} proved the existence and uniqueness of global solution for the 3D anisotropic Boussinesq equations with horizontal viscosity and diffusion in the axisymmetric without swirl, based on a losing estimate with velocity fields of Log-Lipschitz continuity and the algebraic identity $
\frac{u^r}r=\partial_z\Delta^{-1}\frac{\omega^\theta}r-2\frac{\partial_r}{r}\Delta^{-1}\partial_z\Delta^{-1}\frac{\omega^\theta}r$.
Danchin-Paicu \cite{9,10} considered the well-posedness of \eqref{eq1} with partial viscosity $\mu>0, \nu=0$ for sufficiently regular
initial data $(u_0,\rho_0)\in\dot{B}_{p,1}^{-1+\frac3p}
\times \dot{B}_{3,1}^0$ with $p\geq3$, and obtained global well-posedness with the critical quantity $\|u_0\|_{L^{3,\infty}}+\mu^{-1}\|\rho_0\|_{L^1}$ is sufficiently small. As for general case with swirl, some results on global existence for the system \eqref{eq1}, see \cite{3,12,16}.

In fact, the system \eqref{eq1} is the special form of the following general Boussinesq system
\begin{eqnarray}\label{eq2}
\begin{cases}
\partial_tu+u\cdot\nabla u-div(\mu(\rho)\nabla u)+\nabla \pi=f(\rho)\\
\partial_t\rho+u\cdot\nabla\rho-div(\nu(\rho)\nabla \rho)=0\\
div~u=0
\end{cases}
\end{eqnarray}
with one buoyancy force function $f(\rho)$, the general viscosity $\mu(\rho)$ and thermal diffusivity $\nu(\rho)$ depending on the temperature $\rho$, which is the most useful model in geophysical fluid dynamics \cite{23}. Lorca-Boldr \cite{18,19} prove the existence of global weak solution and gained the local well-posedness for general data and the global existence of strong solution for small data. For the partial regularity, or the global regularity with respect to small initial data to the equations \eqref{eq2}, see \cite{11,27}.

The purpose of this paper is to study the global strong solutions to the Boussinesq equations \eqref{eq1} based on the work of Liu-Zhang \cite{17}, where they prove the global existence of strong solutions to three-dimensional
incompressible Navier-Stokes equations with solenoidal initial data and these solutions can be expressed as the Fourier series. But, the method does not directly apply to the Boussinesq equations \eqref{eq1}, owing to buoyancy force function $\rho e_3$. When deal with the pressure term $\pi(t,x)$ and velocity $u(t,x)$, it will lead to estimate $\|r^2f\|_{L_t^1(L^p(r^{1-\frac3p}dx))}$ and $\|r^2f\|_{L_t^{\frac{10}{7-p}}(L^p(r^{1-\frac3p}dx))}$ with $p\in(5,6), r=\sqrt{x_1^2+x_2^2}$ for some function $f(t,x)$ depending on $\rho(t,x)$. It is very difficult to control. In order to overcome the question, we consider the special temperature $\rho(t,x)=(1+r^2)^{-1}\eta(t,x)$. As a result of that the system \eqref{eq1} is modified as a new equations \eqref{eq4} what is more complex than the equations \eqref{eq1} and does not keep some well results of origin, in particular the $L^2$ energy estimation. Compared with \cite{17}, these are some different skills to solve the equations \eqref{eq4}, due to new nonlinear terms. More detail, see Section 4 and 5.

In the end of section, we simple describe the organization of the paper. In Section 2, we introduce some nations and state our main results. In Section 3, we outline our general strategy to prove nonlinear stability. In Section 4 and 5, we give the proofs of the key lemma and proposition.

\section{The main results}
Let us consider the following solutions to the problem \eqref{eq1} in the three-dimensional cylindrical coordinates $(r,\theta, z)$:
\begin{eqnarray}\label{eq3}
u(t,x)=u^r(t,r,\theta,z)e_r+u^\theta(t,r,\theta,z)e_\theta+u^z(t,r,\theta,z)e_z,~~\pi(t,x)=\pi(t,r,\theta,z),~~\rho(t,x)=\frac{\eta(t,r,\theta,z)}{1+r^2},
\end{eqnarray}
where $e_r=(\cos\theta,\sin\theta,0), e_\theta=(-\sin\theta,\cos\theta,0)$ and $e_z=(0,0,1)$. In this case, according to the system \eqref{eq1}, the unknown functions $(u^r,u^\theta,u^z,\pi,\eta)(t,r,\theta,z)$ satisfies the following Cauchy problem
\begin{eqnarray}\label{eq4}
\begin{cases}
D_tu^r-(\partial_r^2+\frac1r\partial_r+\frac1{r^2}\partial_\theta^2+\partial_z^2-\frac1{r^2})u^r-\frac1r(u^\theta)^2+\frac2{r^2}\partial_\theta u^\theta+\partial_r\pi=0,~~~&(t,r,\theta,z)\in(0,T]\times\Omega,\\
D_tu^\theta-(\partial_r^2+\frac1r\partial_r+\frac1{r^2}\partial_\theta^2+\partial_z^2-\frac1{r^2})u^\theta+\frac1ru^ru^\theta-\frac2{r^2}\partial_\theta u^r+\frac1r\partial_\theta\pi=0,~~~&(t,r,\theta,z)\in(0,T]\times\Omega,\\
D_tu^z-(\partial_r^2+\frac1r\partial_r+\frac1{r^2}\partial_\theta^2+\partial_z^2) u^z+\partial_z\pi=\frac{\eta}{1+r^2},~~~&(t,r,\theta,z)\in(0,T]\times\Omega,\\
D_t\eta-(\partial_r^2+\frac1r\partial_r+\frac1{r^2}\partial_\theta^2+\partial_z^2)\eta-\frac{2r}{1+r^2}u^r\eta+\frac{4 r}{1+r^2}\partial_r\eta+\frac{4(1-r^2)}{(1+r^2)^2}\eta=0,~~~&(t,r,\theta,z)\in(0,T]\times\Omega,\\
\partial_ru^r+\frac1ru^r+\frac1r\partial_\theta u^\theta+\partial_zu^z=0,~~~&(t,r,\theta,z)\in(0,T]\times\Omega,\\
(u^r,u^\theta,u^z,\eta)(0,r,\theta,z)=(u_0^r,u_0^\theta,u_0^z,\eta_0)(r,\theta,z),~~~&(r,\theta,z)\in\Omega,
\end{cases}
\end{eqnarray}
where the operator $D_t=\partial_t+u^r\partial_r+\frac1ru^\theta\partial_\theta+u^z\partial_z$ and $\Omega=\mathbb{R}^+\times[0,2\pi]\times \mathbb{R}$.

So as to state our results, we introduce the following notations and the equations.
Denote
\begin{eqnarray*}
\mathcal{M}:=\{f(r,z):\|f\|_{\mathcal{M}}=\|r^{\frac12}f\|_{L^6}+\|f\|_{L^2}+\|(\partial_rf,\partial_zf)\|_{L^2}+\|r^{-1}f\|_{L^2}<+\infty\}.
\end{eqnarray*}
Here and in the following, the norm $\|\cdot\|_{L^p}$ denotes that of $L^p(\mathbb{R}^3)$ for $p>1$. And denote
\begin{align}\label{eq5}
\tilde{U}&=U^re_r+U^ze_z,~~~\tilde{U}_k=U_k^re_r+U_k^ze_z,~~~\tilde{V}_k=V_k^re_r+V_k^ze_z,~~~k\in\mathbb{N}^+,\notag\\
U&=U^re_r+U^\theta e_\theta+U^ze_z,~~~U_k=U_k^re_r+U_k^\theta e_\theta+U_k^ze_z,~~~V_k=V_k^re_r+V_k^\theta e_\theta+V_k^ze_z,~~~k\in\mathbb{N}^+\notag\\
\tilde{\nabla}&=e_r\partial_r+e_z\partial_z,~~~\tilde{\Delta}=\partial_r^2+\frac1r\partial_r+\partial_z^2,~~~A_t=\partial_t+U^r\partial_r+U^z\partial_z,\\
L_t^h&=A_t-\tilde{\Delta}+\frac1{r^2},~~L_t^v=A_t-\tilde{\Delta},~~~\tilde{L}=\frac{4 r}{1+r^2}\partial_r+\frac{4(1-r^2)}{(1+r^2)^2}.\notag
\end{align}
Let the function $(U^r, U^\theta,U^z,\Pi, \xi)$ be one solution to the following system
\begin{eqnarray}\label{eq6}
\begin{cases}
A_tU^r-(\tilde{\Delta}-\frac1{r^2})U^r-\frac1r(U^\theta)^2+\partial_r\Pi=-H^r,~~~&(t,r,z)\in(0,T]\times\Sigma,\\
A_tU^\theta-(\tilde{\Delta}-\frac1{r^2})U^\theta+\frac1rU^rU^\theta=-H^\theta,~~~&(t,r,z)\in(0,T]\times\Sigma,\\
A_tU^z-\tilde{\Delta}U^z+\partial_z\Pi=\frac1{1+r^2}\xi-H^z,~~~&(t,r,z)\in(0,T]\times\Sigma,\\
A_t\xi-(\tilde{\Delta}+\partial_r^2)\xi-\frac{2r}{1+r^2}U^r\xi+\frac{4 r}{1+r^2}\partial_r\xi+\frac{4(1-r^2)}{(1+r^2)^2}\xi=-\phi,~~~&(t,r,z)\in(0,T]\times\Sigma,\\
\partial_rU^r+\frac1rU^r+\partial_zU^z=0,~~~&(t,r,z)\in(0,T]\times\Sigma,\\
(U^r,U^\theta,U^z,\xi)(0,r,z)=(0,0,0),~~~&(r,z)\in\Sigma,
\end{cases}
\end{eqnarray}
where $\Sigma=\mathbb{R}^+\times\mathbb{R}$ and
\begin{align}\label{eq7}
H^r&=\frac12\sum_{k=1}^\infty(\tilde{U}_k\cdot\tilde{\nabla}U_k^r+\tilde{V}_k\cdot\tilde{\nabla}V_k^r-\frac{kN}rU_k^\theta V_k^r+\frac{kN}rV_k^\theta U_k^r)-\frac12\sum_{k=1}^\infty(\frac1r(U_k^\theta)^2+\frac1r(V_k^\theta)^2),\notag\\
H^\theta&=\frac12\sum_{k=1}^\infty(\tilde{U}_k\cdot\tilde{\nabla}U_k^\theta+\tilde{V}_k\cdot\tilde{\nabla}V_k^\theta)+\frac12\sum_{k=1}^\infty(\frac1rU_k^rU_k^\theta+\frac1rV_k^rV_k^\theta),\\
H^z&=\frac12\sum_{k=1}^\infty(\tilde{U}_k\cdot\tilde{\nabla}U_k^z+\tilde{V}_k\cdot\tilde{\nabla}V_k^z-\frac{kN}rU_k^\theta V_k^z+\frac{kN}rV_k^\theta U_k^z),\notag\\
\phi&=\frac12\sum_{k=1}^\infty(\tilde{U}_k\cdot\tilde{\nabla}\xi_k+\tilde{V}_k\cdot\tilde{\nabla}\zeta_k-\frac{kN}rU_k^\theta \zeta_k+\frac{kN}rV_k^\theta \xi_k)-\frac12\sum_{k=1}^\infty(\frac{2r}{1+r^2}U_k^r\xi_k+\frac{2r}{1+r^2}V_k^r\zeta_k).\notag
\end{align}
For $k\in \mathbb{N}^+$, let the function $(U_k^r,U_k^\theta,U_k^z,V_k^r,V_k^\theta,V_k^z,Q_k,R_k, \xi_k,\zeta_k)$ be one solution to the following system
\begin{eqnarray}\label{eq8}
\begin{cases}
L_t^hU_k^r+\tilde{U}_k\cdot\tilde{\nabla}U^r-\frac{kN}rU^\theta V_k^r+\frac{(kN)^2}{r^2}U_k^r-\frac2rU_k^\theta U^\theta-\frac{2kN}{r^2}V_k^\theta+\partial_rQ_k=-F_k^r,\\
L_t^hV_k^r+\tilde{V}_k\cdot\tilde{\nabla}U^r+\frac{kN}rU^\theta U_k^r+\frac{(kN)^2}{r^2}V_k^r-\frac2rU^\theta V_k^\theta+\frac{2kN}{r^2}U_k^\theta+\partial_rR_k=-G_k^r,\\
L_t^hU_k^\theta+\tilde{U}_k\cdot\tilde{\nabla}U^\theta-\frac{kN}rU^\theta V_k^\theta+\frac{(kN)^2}{r^2}U_k^\theta+\frac1r(U^rU_k^\theta+U_k^rU^\theta)+\frac{2kN}{r^2}V_k^r-\frac{kN}rR_k=-F_k^\theta,\\
L_t^hV_k^\theta+\tilde{V}_k\cdot\tilde{\nabla}U^\theta+\frac{kN}rU^\theta U_k^\theta+\frac{(kN)^2}{r^2}V_k^\theta+\frac1r(U^rV_k^\theta+V_k^rU^\theta)-\frac{2kN}{r^2}U_k^r+\frac{kN}rQ_k=-G_k^\theta,\\
L_t^vU_k^z+\tilde{U}_k\cdot\tilde{\nabla}U^z-\frac{kN}rU^\theta V_k^z+\frac{(kN)^2}{r^2}U_k^z+\partial_zQ_k=\frac1{1+r^2}\xi_k-F_k^z,\\
L_t^vV_k^z+\tilde{V}_k\cdot\tilde{\nabla}U^z+\frac{kN}rU^\theta U_k^z+\frac{(kN)^2}{r^2}V_k^z+\partial_zR_k=\frac1{1+r^2}\zeta_k-G_k^z,\\
L_t^v\xi_k+\tilde{U}_k\cdot\tilde{\nabla}\xi-\frac{kN}rU^\theta\zeta_k+\frac{(kN)^2}{r^2}\xi_k-\frac{2r}{1+r^2}U^r\xi_k-\frac{2r}{1+r^2}U_k^r\xi+\tilde{L}\xi_k=-\phi_k,\\
L_t^v\zeta_k+\tilde{V}_k\cdot\tilde{\nabla}\xi+\frac{kN}rU^\theta \xi_k+\frac{(kN)^2}{r^2}\zeta_k-\frac{2r}{1+r^2}U^r\zeta_k-\frac{2r}{1+r^2}V_k^r\xi+\tilde{L}\zeta_k=-\psi_k,\\
\partial_rU_k^r+\frac1rU_k^r+\partial_zU_k^z-\frac{kN}rV_k^\theta=0,~~\partial_rV_k^r+\frac1rV_k^r+\partial_zV_k^z+\frac{kN}rU_k^\theta=0,\\
(U_1^r,V_1^r,U_1^\theta,V_1^\theta,U_1^z,V_1^z,\xi_1,\eta_1)(0,r,z)=(a^r,b^r,N^{-1}a^\theta,N^{-1}b^\theta,a^z,b^z,c,d),\\
(U_k^r,V_k^r,U_k^\theta,V_k^\theta,U_k^z,V_k^z,\xi_k,\eta_k)(0,r,z)=(0,0,0,0,0,0,0,0),~k\geq2
\end{cases}
\end{eqnarray}
where $ t\in(0,T], (r,z)\in\Sigma$, the initial data $(a^r,a^\theta,a^z,b^r,b^\theta,b^z,c,d)(r,z)$ is the given function,  and
\begin{align}\label{eq9}
F_k^r&=\frac12\sum_{k_1+k_2=k}\left(\tilde{U}_{k_1}\cdot\tilde{\nabla}V_{k_2}^r+\frac{k_2N}rU_{k_1}^\theta U_{k_2}^r-\frac1rU_{k_1}^\theta V_{k_2}^\theta+\tilde{V}_{k_1}\cdot\tilde{\nabla}U_{k_2}^r-\frac{k_2N}rV_{k_1}^\theta V_{k_2}^r-\frac1rV_{k_1}^\theta U_{k_2}^\theta\right)\notag\\
&~~~+\frac12\sum_{k_1-k_2=k}\left(\tilde{U}_{k_1}\cdot\tilde{\nabla}V_{k_2}^r+\frac{k_2N}rU_{k_1}^\theta U_{k_2}^r-\frac1rU_{k_1}^\theta V_{k_2}^\theta-\tilde{V}_{k_1}\cdot\tilde{\nabla}U_{k_2}^r+\frac{k_2N}rV_{k_1}^\theta V_{k_2}^r+\frac1rV_{k_1}^\theta U_{k_2}^\theta\right)\notag\\
&~~~+\frac12\sum_{k_2-k_1=k}\left(-\tilde{U}_{k_1}\cdot\tilde{\nabla}V_{k_2}^r-\frac{k_2N}rU_{k_1}^\theta U_{k_2}^r+\frac1rU_{k_1}^\theta V_{k_2}^\theta+\tilde{V}_{k_1}\cdot\tilde{\nabla}U_{k_2}^r-\frac{k_2N}rV_{k_1}^\theta V_{k_2}^r-\frac1rV_{k_1}^\theta U_{k_2}^\theta\right),\notag\\
G_k^r&=\frac12\sum_{k_1+k_2=k}\left(-\tilde{U}_{k_1}\cdot\tilde{\nabla}U_{k_2}^r+\frac{k_2N}rU_{k_1}^\theta V_{k_2}^r+\frac1rU_{k_1}^\theta U_{k_2}^\theta+\tilde{V}_{k_1}\cdot\tilde{\nabla}V_{k_2}^r+\frac{k_2N}rV_{k_1}^\theta U_{k_2}^r-\frac1rV_{k_1}^\theta V_{k_2}^\theta\right)\notag\\
&~~~+\frac12\sum_{|k_1-k_2|=k}\left(\tilde{U}_{k_1}\cdot\tilde{\nabla}U_{k_2}^r-\frac{k_2N}rU_{k_1}^\theta V_{k_2}^r-\frac1rU_{k_1}^\theta U_{k_2}^\theta+\tilde{V}_{k_1}\cdot\tilde{\nabla}V_{k_2}^r+\frac{k_2N}rV_{k_1}^\theta U_{k_2}^r-\frac1rV_{k_1}^\theta V_{k_2}^\theta\right),\notag\\
F_k^\theta&=\frac12\sum_{k_1+k_2=k}\left(\tilde{U}_{k_1}\cdot\tilde{\nabla}V_{k_2}^\theta+\frac{k_2N}rU_{k_1}^\theta U_{k_2}^\theta+\frac1rU_{k_1}^r V_{k_2}^\theta+\tilde{V}_{k_1}\cdot\tilde{\nabla}U_{k_2}^\theta-\frac{k_2N}rV_{k_1}^\theta V_{k_2}^\theta+\frac1rV_{k_1}^r U_{k_2}^\theta\right)\notag\\
&~~~+\frac12\sum_{k_1-k_2=k}\left(\tilde{U}_{k_1}\cdot\tilde{\nabla}V_{k_2}^\theta+\frac{k_2N}rU_{k_1}^\theta U_{k_2}^\theta+\frac1rU_{k_1}^r V_{k_2}^\theta-\tilde{V}_{k_1}\cdot\tilde{\nabla}U_{k_2}^\theta+\frac{k_2N}rV_{k_1}^\theta V_{k_2}^\theta-\frac1rV_{k_1}^r U_{k_2}^\theta\right)\notag\\
&~~~+\frac12\sum_{k_2-k_1=k}\left(-\tilde{U}_{k_1}\cdot\tilde{\nabla}V_{k_2}^\theta-\frac{k_2N}rU_{k_1}^\theta U_{k_2}^\theta-\frac1rU_{k_1}^r V_{k_2}^\theta+\tilde{V}_{k_1}\cdot\tilde{\nabla}U_{k_2}^\theta-\frac{k_2N}rV_{k_1}^\theta V_{k_2}^\theta+\frac1rV_{k_1}^r U_{k_2}^\theta\right),\\
G_k^\theta&=\frac12\sum_{k_1+k_2=k}\left(-\tilde{U}_{k_1}\cdot\tilde{\nabla}U_{k_2}^\theta+\frac{k_2N}rU_{k_1}^\theta V_{k_2}^\theta-\frac1rU_{k_1}^r U_{k_2}^\theta+\tilde{V}_{k_1}\cdot\tilde{\nabla}V_{k_2}^\theta+\frac{k_2N}rV_{k_1}^\theta U_{k_2}^\theta+\frac1rV_{k_1}^r V_{k_2}^\theta\right)\notag\\
&~~~+\frac12\sum_{|k_1-k_2|=k}\left(\tilde{U}_{k_1}\cdot\tilde{\nabla}U_{k_2}^\theta-\frac{k_2N}rU_{k_1}^\theta V_{k_2}^\theta+\frac1rU_{k_1}^r U_{k_2}^\theta+\tilde{V}_{k_1}\cdot\tilde{\nabla}V_{k_2}^\theta+\frac{k_2N}rV_{k_1}^\theta U_{k_2}^\theta+\frac1rV_{k_1}^r V_{k_2}^\theta\right),\notag\\
F_k^z&=\frac12\sum_{k_1+k_2=k}\left(\tilde{U}_{k_1}\cdot\tilde{\nabla}V_{k_2}^z+\frac{k_2N}rU_{k_1}^\theta U_{k_2}^z+\tilde{V}_{k_1}\cdot\tilde{\nabla}U_{k_2}^z-\frac{k_2N}rV_{k_1}^\theta V_{k_2}^z\right)\notag\\
&~~~+\frac12\sum_{k_1-k_2=k}\left(\tilde{U}_{k_1}\cdot\tilde{\nabla}V_{k_2}^z+\frac{k_2N}rU_{k_1}^\theta U_{k_2}^z-\tilde{V}_{k_1}\cdot\tilde{\nabla}U_{k_2}^z+\frac{k_2N}rV_{k_1}^\theta V_{k_2}^z\right)\notag\\
&~~~+\frac12\sum_{k_2-k_1=k}\left(-\tilde{U}_{k_1}\cdot\tilde{\nabla}V_{k_2}^z-\frac{k_2N}rU_{k_1}^\theta U_{k_2}^z+\tilde{V}_{k_1}\cdot\tilde{\nabla}U_{k_2}^z-\frac{k_2N}rV_{k_1}^\theta V_{k_2}^z\right),\notag\\
G_k^z&=\frac12\sum_{k_1+k_2=k}\left(-\tilde{U}_{k_1}\cdot\tilde{\nabla}U_{k_2}^z+\frac{k_2N}rU_{k_1}^\theta V_{k_2}^z+\tilde{V}_{k_1}\cdot\tilde{\nabla}V_{k_2}^z+\frac{k_2N}rV_{k_1}^\theta U_{k_2}^z\right)\notag\\
&~~~+\frac12\sum_{|k_1-k_2|=k}\left(\tilde{U}_{k_1}\cdot\tilde{\nabla}U_{k_2}^z-\frac{k_2N}rU_{k_1}^\theta V_{k_2}^z+\tilde{V}_{k_1}\cdot\tilde{\nabla}V_{k_2}^z+\frac{k_2N}rV_{k_1}^\theta U_{k_2}^z\right),\notag\\
\phi_k&=\frac12\sum_{k_1+k_2=k}\left(\tilde{U}_{k_1}\cdot\tilde{\nabla}\zeta_{k_2}+\frac{k_2N}rU_{k_1}^\theta \xi_{k_2}+\tilde{V}_{k_1}\cdot\tilde{\nabla}\xi_{k_2}-\frac{k_2N}rV_{k_1}^\theta \zeta_{k_2}-\frac{2r}{1+r^2}U_{k_1}^r\zeta_{k_2}-\frac{2r}{1+r^2}V_{k_1}^r\xi_{k_2}\right)\notag\\
&~~~+\frac12\sum_{k_1-k_2=k}\left(\tilde{U}_{k_1}\cdot\tilde{\nabla}\zeta_{k_2}+\frac{k_2N}rU_{k_1}^\theta \xi_{k_2}-\tilde{V}_{k_1}\cdot\tilde{\nabla}\xi_{k_2}+\frac{k_2N}rV_{k_1}^\theta \zeta_{k_2}-\frac{2r}{1+r^2}U_{k_1}^r\zeta_{k_2}+\frac{2r}{1+r^2}V_{k_1}^r\xi_{k_2}\right)\notag\\
&~~~+\frac12\sum_{k_2-k_1=k}\left(-\tilde{U}_{k_1}\cdot\tilde{\nabla}\zeta_{k_2}-\frac{k_2N}rU_{k_1}^\theta \xi_{k_2}+\tilde{V}_{k_1}\cdot\tilde{\nabla}\xi_{k_2}-\frac{k_2N}rV_{k_1}^\theta \zeta_{k_2}+\frac{2r}{1+r^2}U_{k_1}^r\zeta_{k_2}-\frac{2r}{1+r^2}V_{k_1}^r\xi_{k_2}\right),\notag\\
\psi_k&=\frac12\sum_{k_1+k_2=k}\left(-\tilde{U}_{k_1}\cdot\tilde{\nabla}\xi_{k_2}+\frac{k_2N}rU_{k_1}^\theta \zeta_{k_2}+\tilde{V}_{k_1}\cdot\tilde{\nabla}\zeta_{k_2}+\frac{k_2N}rV_{k_1}^\theta \xi_{k_2}+\frac{2r}{1+r^2}U_{k_1}^r\xi_{k_2}-\frac{2r}{1+r^2}V_{k_1}^r\zeta_{k_2}\right)\notag\\
&~~~+\frac12\sum_{|k_1-k_2|=k}\left(\tilde{U}_{k_1}\cdot\tilde{\nabla}\xi_{k_2}-\frac{k_2N}rU_{k_1}^\theta \zeta_{k_2}+\tilde{V}_{k_1}\cdot\tilde{\nabla}\zeta_{k_2}+\frac{k_2N}rV_{k_1}^\theta \xi_{k_2}-\frac{2r}{1+r^2}U_{k_1}^r\xi_{k_2}-\frac{2r}{1+r^2}V_{k_1}^r\zeta_{k_2}\right).\notag
\end{align}
For the problem \eqref{eq1}, we have the following result, which is the main result of this paper.

\textbf{Theorem 2.1} For any $(a^r,a^\theta,a^z,b^r,b^\theta,b^z,c,d)\in\mathcal{M}$ satisfying
\begin{eqnarray*}
\partial_ra^r+\frac1ra^r+\partial_za^z+\frac1rb^\theta=0,~~~\partial_rb^r+\frac1rb^r+\partial_zb^z-\frac1ra^\theta=0
\end{eqnarray*}
and any given finite time $T>0$, there exists an integer $N_0>0$ depending only on $\|(a^r,a^\theta,a^z,b^r,b^\theta,b^z,c,d)\|_{\mathcal{M}}$ and $T$ such that for any integer $N\geq N_0$, the system \eqref{eq1} with $\mu=\nu=1$ and initial data
\begin{eqnarray}\label{initial-s}
\begin{cases}
(u_0,\rho_0)(x)=(u ^r_0(r,\theta,z)e_r+u^\theta_0(r,\theta,z)e_\theta+u^z_0(r,\theta,z)e_z, \frac{\eta_0(r,\theta,z)}{1+r^2}),\\
u_0^r(r,\theta,z)=a^r(r,z)\sin(N\theta)+b^r(r,z)\cos(N\theta),\\
u_0^\theta(r,\theta,z)=\frac1{N}(a^\theta(r,z)\sin(N\theta)+b^\theta(r,z)\cos(N\theta)),\\
u_0^z(r,\theta,z)=a^z(r,z)\sin(N\theta)+b^z(r,z)\cos(N\theta),\\
\eta_0(r,\theta,z)=c(r,z)\sin(N\theta)+d(r,z)\cos(N\theta)
\end{cases}
\end{eqnarray}
has a unique global strong solution defined on the interval $[0,T]$, which can be expressed into the Fourier series: for all $t\in[0,T]$
\begin{eqnarray}\label{fourier-s}
\begin{cases}
u^r(t,r,\theta,z)=U^r(t,r,z)+\sum\limits_{k=1}^\infty(U_k^r(t,r,z)\sin(kN\theta)+ V_k^r(t,r,z)\cos(kN\theta)),\\
u^\theta(t,r,\theta,z)=U^\theta(t,r,z)+\sum\limits_{k=1}^\infty(U_k^\theta(t,r,z)\sin(kN\theta)+ V_k^\theta(t,r,z)\cos(kN\theta)),\\
u^z(t,r,\theta,z)=U^z(t,r,z)+\sum\limits_{k=1}^\infty(U_k^z(t,r,z)\sin(kN\theta)+ V_k^z(t,r,z)\cos(kN\theta)),\\
\pi(t,r,\theta,z)=\Pi(t,r,z)+\sum\limits_{k=1}^\infty(Q_k(t,r,z)\sin(kN\theta) +R_k(t,r,z)\cos(kN\theta)),\\
\rho(t,r,\theta,z)=\frac1{1+r^2}\left(\xi(t,r,z)+\sum\limits_{k=1}^\infty(\xi_k(t,r,z)\sin(kN\theta)+ \zeta_k(t,r,z)\cos(kN\theta))\right),
\end{cases}
\end{eqnarray}
where the lead term $(U^r,U^\theta,U^z,\Pi, \xi)$ and the general $k-$order term $(U_k^r,V_k^r,U_k^\theta,V_k^\theta,U_k^z,V_k^z,Q_k,R_k,\xi_k,\zeta_k), k\in N^+,$ satisfy the system \eqref{eq6}-\eqref{eq7} and the system \eqref{eq8}-\eqref{eq9} respectively.

Furthermore, for any $p\in(5,6)$ and $\alpha_p,\beta_p$ satisfying
\begin{eqnarray*}
1<\beta_p<\frac{p-1}4,~~~0<4\alpha_p<p-3-2\beta_p,
\end{eqnarray*}
there exists a positive constant $C_{in}(T)$ depending only $(a^r,a^\theta,a^z,b^r,b^\theta,b^z,c,d)$ and $T$ such that Fourier coefficients in \eqref{fourier-s} satisfy the following estimates:
\begin{align*}
&\|U^re_r+U^\theta e_\theta+U^ze_z\|_{L^\infty([0,T];L^3)}+\|\xi\|_{L^\infty([0,T];L^3)}\leq C_{in}(T)N^{-\frac15},\\
&\|r^{1-\frac3p}\varpi_1\|_{L^\infty([0,T];L^p)}^p\leq C_{in}(T),~~\sum_{k=2}^\infty(kN)^{2\beta_p}\|r^{1-\frac3p}\varpi_k\|_{L^\infty([0,T]; L^p)}^p\leq C_{in}(T)N^{-2\alpha_p},
\end{align*}
where $\varpi_k=(U_k^r,V_k^r,U_k^z,V_k^z,\xi_k,\zeta_k,\sqrt{kN}U_k^\theta,\sqrt{kN}V_k^\theta)$ for all $k\in\mathbb{N}^+$.

\textbf{Remark 2.1:} The result in Theorem 1.1 is one result on the special large strong solution $u$ with small swirl $u^\theta$ because the initial $u_0^\theta$ is small when $N$ is large enough. In fact, if
\begin{eqnarray*}
U_k^r=U_k^z=V_k^\theta=Q_k=\xi_k\equiv0
\end{eqnarray*}
hold for all $k\in N^+$, then
\begin{eqnarray*}
H^\theta=F_k^r=F_k^z=G_k^\theta=\phi_k\equiv0
\end{eqnarray*}
for any $k\in N^+$. Thus, we can take $U^\theta\equiv0$ in the Fourier series \eqref{fourier-s} if one assume $a^r=a^z=a^\theta=b^\theta=c=0$, and $U^\theta$ is also small enough in the case of general initial data when $N$ is large enough. Moreover, the swirls $U_1^\theta$ and $V_1^\theta$ in the Fourier series \eqref{fourier-s}  have the structure:
\begin{eqnarray*}
U_1^\theta(t,r,z)=\frac1N\bar{U}_1^\theta(t,r,z),~~~V_1^\theta(t,r,z)=\frac1N\bar{V}_1^\theta(t,r,z) \end{eqnarray*}
with
\begin{eqnarray*}
\|r^{1-\frac3p}(\bar{U}_1^\theta,\bar{V}_1^\theta)\|_{L^\infty([0,T];L^p)}\leq C_{in}(T)N^{\frac12},
\end{eqnarray*} and the high order Fourier coefficients $U_k^\theta$ and $V_k^\theta$, $k\ge 2$, in the Fourier series \eqref{fourier-s} have similar smallness structure.

\textbf{Remark 2.2:} If the continuous function $g(r)$ satisfies that for all $r\in\mathbb{R}^+$ and $\alpha=1,2$,
\begin{eqnarray*}
0<(1+r^2)g(r)<+\infty,~~~|\frac{r^\alpha\partial_r^\alpha g(r)}{g(r)}|<+\infty
\end{eqnarray*}
then, the equations \eqref{eq1} has the general solutions:
\begin{eqnarray*}
u(t,x)=u^r(t,r,\theta,z)e_r+u^\theta(t,r,\theta,z)e_\theta+u^z(t,r,\theta,z)e_z,~~\pi(t,x)=\Pi(t, r,\theta,z),~~\rho(t,x)=g(r)\eta(t,r,\theta,z),
\end{eqnarray*}
which can be expressed as the similar Fourier series to \eqref{fourier-s}. Thus, we can extend our result in Theorem 2.1 to the system \eqref{eq2} with $f(\rho)=g(r)\rho e_3$ and $\mu(\rho)=\nu(\rho)=1$.

\section{The proof of Theorem 2.1}
In this section, we are devoted to proving Theorem 2.1. Based on the structure \eqref{eq3} of the solution to the system \eqref{eq1}, we just only need to obtain the global well-posedness of the system \eqref{eq4} with suitable smooth initial data \eqref{eq12}. Motivated by the work \cite{17}, for any finite time $T>0$, we try to consider the following Fourier series solutions to the system \eqref{eq4} satisfying the initial data \eqref{eq12}: for all $t\in(0,T]$,
\begin{eqnarray}\label{eq10}
\begin{cases}
u^r(t,r,\theta,z)=U^r(t,r,z)+\sum\limits_{k=1}^\infty(U_k^r(t,r,z)\sin(kN\theta)+ V_k^r(t,r,z)\cos(kN\theta)),\\
u^\theta(t,r,\theta,z)=U^\theta(t,r,z)+\sum\limits_{k=1}^\infty(U_k^\theta(t,r,z)\sin(kN\theta)+ V_k^\theta(t,r,z)\cos(kN\theta)),\\
u^z(t,r,\theta,z)=U^z(t,r,z)+\sum\limits_{k=1}^\infty(U_k^z(t,r,z)\sin(kN\theta)+ V_k^z(t,r,z)\cos(kN\theta)),\\
\eta(t,r,\theta,z)=\xi(t,r,z)+\sum\limits_{k=1}^\infty(\xi_k(t,r,z)\sin(kN\theta)+ \zeta_k(t,r,z)\cos(kN\theta)),
\end{cases}
\end{eqnarray}
and
\begin{eqnarray}\label{eq11}
\pi(t,r,\theta,z)=\Pi(t,r,z)+\sum_{k=1}^\infty(Q_k(t,r,z)\sin(kN\theta) +R_k(t,r,z)\cos(kN\theta)).
\end{eqnarray}
Once it can achieve, our result can be justified right away. In fact, for any the local smooth solution of the system \eqref{eq4}, the above expressions are right.

Using the facts that for all $k_1,k_2\in\mathbb{N}^+$
\begin{eqnarray*}
2\cos(k_1\theta)\cdot\cos(k_2\theta)=\cos(k_1+k_2)\theta+\cos(k_1-k_2)\theta,\\
2\sin(k_1\theta)\cdot\sin(k_2\theta)=\cos(k_1-k_2)\theta-\cos(k_1+k_2)\theta,\\
2\sin(k_1\theta)\cdot\cos(k_2\theta)=\sin(k_1+k_2)\theta+\sin(k_1-k_2)\theta
\end{eqnarray*}
and submitting the solutions \eqref{eq10} and \eqref{eq11} into \eqref{eq4} with initial data \eqref{initial-s}, we obtain the system \eqref{eq6}-\eqref{eq7} for the lead term $(U^r,U^\theta,U^z,\Pi, \xi)$ in the Fourier series \eqref{fourier-s} and the system \eqref{eq8}-\eqref{eq9} the general $k-$order term $(U_k^r,\cdots,\zeta_k), k\in N^+,$ in the Fourier series \eqref{fourier-s} respectively.

In the rest sections, we will prove that the system \eqref{eq4} has the solution \eqref{eq10}-\eqref{eq11} with initial data
\begin{eqnarray}\label{eq12}
\begin{cases}
u_0^r(r,\theta,z)=a^r(r,z)\sin(N\theta)+b^r(r,z)\cos(N\theta),\\
u_0^\theta(r,\theta,z)=\frac1{N}(a^\theta(r,z)\sin(N\theta)+b^\theta(r,z)\cos(N\theta)),\\
u_0^z(r,\theta,z)=a^z(r,z)\sin(N\theta)+b^z(r,z)\cos(N\theta),\\
\eta_0(r,\theta,z)=c(r,z)\sin(N\theta)+d(r,z)\cos(N\theta)
\end{cases}
\end{eqnarray}
satisfying some a priori estimates, which are given in Proposition 3.1 and Proposition 3.2 below. In order to state them, we introduce the foloowing energy functionals:
\begin{align*}
E_p(t)&=N^{-2\alpha_p}\left(\|r^{\frac{p-3}2}|\tilde{U}_1(t)|^{\frac{p}2}\|_{L^2}^2+\int_0^t\|r^{\frac{p-3}2}|\tilde{\nabla}\tilde{U}_1||\tilde{U}_1|^{\frac{p}2-1}\|_{L^2}^2
+N^2\int_0^t\|r^{\frac{p-5}2}|\tilde{U}_1|^{\frac{p}2}\|_{L^2}^2dt'\right)\\
&+\sum_{k=2}^\infty(kN)^{2\beta_p}\left(\|r^{\frac{p-3}2}|\tilde{U}_k(t)|^{\frac{p}2}\|_{L^2}^2+\int_0^t\|r^{\frac{p-3}2}|\tilde{\nabla}\tilde{U}_k||\tilde{U}_k|^{\frac{p}2-1}\|_{L^2}^2+(kN)^2\int_0^t\|r^{\frac{p-5}2}|\tilde{U}_k|^{\frac{p}2}\|_{L^2}^2dt'\right)\\
&+N^{-2\alpha_p}\left(\|r^{\frac{p-3}2}|\tilde{V}_1(t)|^{\frac{p}2}\|_{L^2}^2+\int_0^t\|r^{\frac{p-3}2}|\tilde{\nabla}\tilde{V}_1||\tilde{V}_1|^{\frac{p}2-1}\|_{L^2}^2
+N^2\int_0^t\|r^{\frac{p-5}2}|\tilde{V}_1|^{\frac{p}2}\|_{L^2}^2dt'\right)\\
&+\sum_{k=2}^\infty(kN)^{2\beta_p}\left(\|r^{\frac{p-3}2}|\tilde{V}_k(t)|^{\frac{p}2}\|_{L^2}^2+\int_0^t\|r^{\frac{p-3}2}|\tilde{\nabla}\tilde{V}_k||\tilde{V}_k|^{\frac{p}2-1}\|_{L^2}^2+(kN)^2\int_0^t\|r^{\frac{p-5}2}|\tilde{V}_k|^{\frac{p}2}\|_{L^2}^2dt'\right)\\
&+N^{\frac{p}2-2\alpha_p}\left(\|r^{\frac{p-3}2}|U_1^\theta(t)|^{\frac{p}2}\|_{L^2}^2+\int_0^t\|r^{\frac{p-3}2}|\tilde{\nabla}U_1^\theta||U_1^\theta|^{\frac{p}2-1}\|_{L^2}^2
+N^2\int_0^t\|r^{\frac{p-5}2}|U_1^\theta|^{\frac{p}2}\|_{L^2}^2dt'\right)\\
&+\sum_{k=2}^\infty(kN)^{\frac{p}2+2\beta_p}\left(\|r^{\frac{p-3}2}|U_k^\theta(t)|^{\frac{p}2}\|_{L^2}^2+\int_0^t\|r^{\frac{p-3}2}|\tilde{\nabla}U_k^\theta||U_k^\theta|^{\frac{p}2-1}\|_{L^2}^2+(kN)^2\int_0^t\|r^{\frac{p-5}2}|U_k^\theta|^{\frac{p}2}\|_{L^2}^2dt'\right)\\
&+N^{\frac{p}2-2\alpha_p}\left(\|r^{\frac{p-3}2}|V_1^\theta(t)|^{\frac{p}2}\|_{L^2}^2+\int_0^t\|r^{\frac{p-3}2}|\tilde{\nabla}V_1^\theta||V_1^\theta|^{\frac{p}2-1}\|_{L^2}^2
+N^2\int_0^t\|r^{\frac{p-5}2}|V_1^\theta|^{\frac{p}2}\|_{L^2}^2dt'\right)\\
&+\sum_{k=2}^\infty(kN)^{\frac{p}2+2\beta_p}\left(\|r^{\frac{p-3}2}|V_k^\theta(t)|^{\frac{p}2}\|_{L^2}^2+\int_0^t\|r^{\frac{p-3}2}|\tilde{\nabla}V_k^\theta||U_k^\theta|^{\frac{p}2-1}\|_{L^2}^2+(kN)^2\int_0^t\|r^{\frac{p-5}2}|V_k^\theta|^{\frac{p}2}\|_{L^2}^2dt'\right),\\
\mathcal{E}_p(t)&=N^{-2\alpha_p}\left(\|r^{\frac{p-3}2}|\xi_1(t)|^{\frac{p}2}\|_{L^2}^2+\int_0^t\|r^{\frac{p-3}2}|\tilde{\nabla}\xi_1||\xi_1|^{\frac{p}2-1}\|_{L^2}^2
+N^2\int_0^t\|r^{\frac{p-5}2}|\xi_1|^{\frac{p}2}\|_{L^2}^2dt'\right)\\
&+\sum_{k=2}^\infty(kN)^{2\beta_p}\left(\|r^{\frac{p-3}2}|\xi_k(t)|^{\frac{p}2}\|_{L^2}^2+\int_0^t\|r^{\frac{p-3}2}|\tilde{\nabla}\xi_k||\xi_k|^{\frac{p}2-1}\|_{L^2}^2+(kN)^2\int_0^t\|r^{\frac{p-5}2}|\xi_k|^{\frac{p}2}\|_{L^2}^2dt'\right)\\
&+N^{-2\alpha_p}\left(\|r^{\frac{p-3}2}|\zeta_1(t)|^{\frac{p}2}\|_{L^2}^2+\int_0^t\|r^{\frac{p-3}2}|\tilde{\nabla}\zeta_1||\zeta_1|^{\frac{p}2-1}\|_{L^2}^2
+N^2\int_0^t\|r^{\frac{p-5}2}|\zeta_1|^{\frac{p}2}\|_{L^2}^2dt'\right)\\
&+\sum_{k=2}^\infty(kN)^{2\beta_p}\left(\|r^{\frac{p-3}2}|\zeta_k(t)|^{\frac{p}2}\|_{L^2}^2+\int_0^t\|r^{\frac{p-3}2}|\tilde{\nabla}\zeta_k||\zeta_k|^{\frac{p}2-1}\|_{L^2}^2+(kN)^2\int_0^t\|r^{\frac{p-5}2}|\zeta_k|^{\frac{p}2}\|_{L^2}^2dt'\right),
\end{align*}
and
\begin{eqnarray*}
D(t)=\sum_{k=1}^\infty\int_0^t\|r^{-\frac32}U_k\|_{L^2}^2+\|r^{-\frac32}V_k\|_{L^2}^2dt'. \end{eqnarray*}
In Section 4, we shall prove the following a priori estimate.

\textbf{Proposition 3.1:} Let $0<T<+\infty$ and $0<\alpha_p<\frac{p-3-2\beta_p}4$ and $1<\beta_p<\frac{p-1}4$ with $p\in(5,6)$. Assume $(U^r,U^\theta,U^z,\Pi, \xi)$ and $(U_k^r,V_k^r,U_k^\theta,V_k^\theta,U_k^z,V_k^z,Q_k,R_k,\xi_k,\zeta_k)$ be
smooth solutions of the system \eqref{eq6} and the system \eqref{eq8}, defined on $[0,T]$, respectively. Then there exists a positive integer $N_0$ such that for all $N>N_0$ and $0\leq t\leq T$,
\begin{align*}
(\frac12-C(\|U\|_{L_t^\infty(L^3)}+\|\xi\|_{L_t^\infty(L^3)}))(E_p(t)+\mathcal{E}_p(t))&\leq E_p(0)+C(E_p(t)+\mathcal{E}_p(t))^{1+\frac1p}+CN^{-2}D^{\frac{p+1}2}(t)\\
&+CN^{-p}(\|U\|_{L_t^\infty(L^3)}+N^{-\frac{p}{2(p-5)}})D^{\frac{p}2}(t) \end{align*}
holds for some positive constant $C$ independent of the time $T,N$ and the solutions $(U^r,\cdots,\zeta_k), k\ge 1$.

With Proposition 3.1 at hand, we shall establish the a priori estimates of the solution $(u,\eta)$ to the system \eqref{eq4} through the estimates of its Fourier coefficients by using Plancherel's identity and
Hausdorff-Young's inequality. Finally we close the global a priori estimates by a standard continuity argument. Precisely, we shall prove the following proposition in Section 5:

\textbf{Proposition 3.2.} There exists a positive integer $N_0$ to be large enough so that for $N>N_0$ and any finite time $T>0$, the system \eqref{eq4} with initial data \eqref{eq12} admits a unique global strong solution $(u^r, u^\theta, u^z, \Pi, \eta)(t,r,\theta, z)$, defined on $[0,T]\times \mathbb{R}^+\times [0,2\pi]\times \mathbb{R}$ and expressed by \eqref{eq10}-\eqref{eq11}, satisfying the following estimations:
\begin{description}
  \item[i:] component estimates:
  \begin{align*}
  &\|U^re_r+U^\theta e_\theta+U^ze_z\|_{L^\infty([0,T];L^3)}+\|\xi\|_{L^\infty([0,T];L^3)}\leq C_{in}(T)N^{-\frac15},\\
  &\|r^{1-\frac3p}\varpi_1\|_{L^\infty([0,T];L^p)}^p\leq C_{in}(T),~~\sum_{k=2}^\infty(kN)^{2\beta_p}\|r^{1-\frac3p}\varpi_k\|_{L^\infty([0,T]; L^p)}^p\leq C_{in}(T)N^{-2\alpha_p};
  \end{align*}
  \item[ii:] entirety estimates:
  $(u,\eta)\in (L^\infty([0,T],H^1))^2$, $(\nabla u,\nabla\eta)\in (L^2([0,T],H^1))^2$ and
  \begin{eqnarray*}
  \|\nabla u\|_{L_T^\infty(L^2)}^2+\|\nabla \eta\|_{L_T^\infty(L^2)}^2+\|\Delta u\|_{L_T^2(L^2)}^2+\|\Delta \eta\|_{L_T^2(L^2)}^2\leq C_{in}(T),~~~\|u\|_{L_T^5(L^5)}\leq C_{in}(T)N^{-\frac15};
  \end{eqnarray*}
\end{description}
where $\varpi_k=(U_k^r,V_k^r,U_k^z,V_k^z,\xi_k,\zeta_k,\sqrt{kN}U_k^\theta,\sqrt{kN}V_k^\theta)$ for all $k\in\mathbb{N}^+$, and $C_{in}(T)$ is some positive constant depending only on $\|(a^r,a^\theta,a^z,b^r,b^\theta,b^z,c,d)\|_{\mathcal{M}}$ and $T$.

Now, we sketch the proof of Theorem 2.1.

\textbf{Proof of Theorem 2.1:} By Proposition 3.2, we know that the system \eqref{eq4} has one unique global strong solution \eqref{eq10}-\eqref{eq11} of Fourier series type when the initial data are given by \eqref{eq12}. This, together with \eqref{eq3} and the uniqueness of the strong solution to the Cauchy problem to three-dimensional incompressible Boussinesq equations with the smooth initial data, completes the proof of the existence part in Theorem
2.1. The inequalities of Theorem 2.1 is ensured by the component estimates in the Proposition 3.2. This completes the proof of Theorem 2.1. $\square$

\section{The proof of Proposition 3.1}
This section is devoted to the proof of Proposition 3.1, which is the core to the proof of Theorem 2.1. In order to prove the result, we need the following four consequences. The first and second comes from the work \cite{21}, and the rest will be justified at the end of the section.

\textbf{Proposition 4.1}(\cite{21}) Let $m\geq3, 1<q<p<\infty$, and $\alpha,\beta$ satisfy
\begin{eqnarray*}
\alpha+\beta>0,~~~\alpha+\frac1p\geq-3,~~~\beta-\frac1q\geq-2,~~~and~~~\frac1q=\frac1p+\frac{\alpha+\beta}3.
\end{eqnarray*}
Assume $g(r,z)$ is a solution of the equation:
\begin{eqnarray*}
\mathcal{L}_mg:=-(\partial_r^2+\frac1r\partial_r+\partial_z^2)g+\frac{m^2}{r^2}g=f(r,z). \end{eqnarray*}
Then for any axisymmetric function $f(r,z)$ and any sufficiently small $\epsilon>0$, one has
\begin{eqnarray*}
\|r^\alpha \mathcal{L}_m^{-1}(r^{\beta-2}f)\|_{L^p}\leq C(\epsilon) m^{-(1+\frac1p-\frac1q-\epsilon)}\|f\|_{L^q},
\end{eqnarray*}
where the positive number constant $C(\epsilon)$ is a only dependent of $\epsilon$.

\textbf{Lemma 4.1}(\cite{21}) Let $0<\alpha_p<\frac{p-3-2\beta_p}4$ and $1<\beta_p<\frac{p-1}4$ with $p\in(5,6)$. Assume the function $\tilde{Q}_k$ and $\tilde{R}_k$ respectively is a solution to the equations:
\begin{align*}
\mathcal{L}_{kN}\tilde{Q}_k&=(\partial_r+\frac1r)(\tilde{U}\cdot\tilde{\nabla}U_k^r+\tilde{U}_k\cdot\tilde{\nabla}U^r-\frac{kN}rU^\theta V_k^r-\frac2rU_k^\theta U^\theta)\\
&-\frac{kN}r(\tilde{U}\cdot\tilde{\nabla}V_k^\theta+\tilde{V}_k\cdot\tilde{\nabla}U^\theta+\frac{kN}rU^\theta U_k^\theta+\frac1r(U^rV_k^\theta+V_k^rU^\theta))\\
&+\partial_z(\tilde{U}\cdot\tilde{\nabla}U_k^z+\tilde{U}_k\cdot\tilde{\nabla}U^z-\frac{kN}rU^\theta V_k^z)+(\partial_r+\frac1r)F_k^r-\frac{kN}rG_k^\theta+\partial_zF_k^z\\
&-(\partial_r+\frac1r)\tilde{\Delta}U_k^r+\frac{kN}r\tilde{\Delta}V_k^\theta-\partial_z\tilde{\Delta}U_k^z-\frac2{r^3}U_k^r-2kN(\partial_r+\frac1r)(\frac1{r^2}V_k^\theta)-\frac1{r^2}\partial_zU_k^z,
\end{align*}
and
\begin{align*}
\mathcal{L}_{kN}\tilde{R}_k&=(\partial_r+\frac1r)(\tilde{U}\cdot\tilde{\nabla}V_k^r+\tilde{V}_k\cdot\tilde{\nabla}U^r+\frac{kN}rU^\theta U_k^r-\frac2rV_k^\theta U^\theta)\\
&+\frac{kN}r(\tilde{U}\cdot\tilde{\nabla}U_k^\theta+\tilde{U}_k\cdot\tilde{\nabla}U^\theta-\frac{kN}rU^\theta V_k^\theta+\frac1r(U^rU_k^\theta+U_k^rU^\theta))\\
&+\partial_z(\tilde{U}\cdot\tilde{\nabla}V_k^z+\tilde{V}_k\cdot\tilde{\nabla}U^z+\frac{kN}rU^\theta U_k^z)+(\partial_r+\frac1r)G_k^r+\frac{kN}rF_k^\theta+\partial_zG_k^z\\
&-(\partial_r+\frac1r)\tilde{\Delta}V_k^r-\frac{kN}r\tilde{\Delta}U_k^\theta-\partial_z\tilde{\Delta}V_k^z-\frac2{r^3}V_k^r+2kN(\partial_r+\frac1r)(\frac1{r^2}U_k^\theta)-\frac1{r^2}\partial_zV_k^z,
\end{align*}
then
\begin{align*}
&N^{-2\alpha_p}|\int_0^t\int_{\mathbb{R}^3}(\partial_r\tilde{Q}_1U_1^r+\partial_z\tilde{Q}_1U_1^z)|\tilde{U}_1|^{p-2}r^{p-3}+\frac{N}r\tilde{Q}_1 N^{\frac{p}2}V_1^\theta|V_1^\theta|^{p-2}r^{p-3}dxdt'|\\
&+N^{-2\alpha_p}|\int_0^t\int_{\mathbb{R}^3}(\partial_r\tilde{R}V_1^r+\partial_z\tilde{R}V_1^z)|\tilde{V}_1|^{p-2}r^{p-3}-\frac{N}rR_1 N^{\frac{p}2}U_1^\theta|U_1^\theta|^{p-2}r^{p-3}dxdt'|\\
&+\sum_{k\geq2}(kN)^{2\beta_p}|\int_0^t\int_{\mathbb{R}^3}(\partial_r\tilde{Q}_kU_k^r+\partial_z\tilde{Q}_kU_k^z)|\tilde{U}_k|^{p-2}r^{p-3}+\frac{kN}r\tilde{Q}_k (kN)^{\frac{p}2}V_k^\theta|V_k^\theta|^{p-2}r^{p-3}dxdt'|\\
&+\sum_{k\geq2}(kN)^{2\beta_p}|\int_0^t\int_{\mathbb{R}^3}(\partial_r\tilde{R}_kV_k^r+\partial_z\tilde{R}_kV_k^z)|\tilde{V}_k|^{p-2}r^{p-3}-\frac{kN}r\tilde{R}_k (kN)^{\frac{p}2}U_k^\theta|U_k^\theta|^{p-2}r^{p-3}dxdt'|\\
&\leq (\frac14+C\|U\|_{L_t^\infty(L^3)})E_p(t)+CE_p^{1+\frac1p}(t)+CN^{-2}D^{\frac{p+1}2}(t)
+CN^{-p}(\|U\|_{L_t^\infty(L^3)}+N^{-\frac{p}{2(p-5)}})D^{\frac{p}2}(t). \end{align*}

Thanks to the Proposition 4.1, we have the next Lemma 4.2.

\textbf{Lemma 4.2} Let $0<\alpha_p<\frac{p-3-2\beta_p}4$ and $1<\beta_p<\frac{p-1}4$ with $p\in(5,6)$. Assume the function $\bar{Q}_k$ and $\bar{R}_k$ respectively is a solution to the equations:
\begin{eqnarray*}
\mathcal{L}_{kN}\bar{Q}_k=-\frac1{1+r^2}\xi_k,~~~\mathcal{L}_{kN}\bar{R}_k=-\frac1{1+r^2}\zeta_k,
\end{eqnarray*}
then
\begin{align*}
N^{-2\alpha_p}|\int_0^t\int_{\mathbb{R}^3}(\partial_r\bar{Q}_1U_1^r&+\partial_z\bar{Q}_1U_1^z)|\tilde{U}_1|^{p-2}r^{p-3}+\frac{N}r\bar{Q}_1 N^{\frac{p}2}V_1^\theta|V_1^\theta|^{p-2}r^{p-3}dxdt'|\\
&+N^{-2\alpha_p}|\int_0^t\int_{\mathbb{R}^3}(\partial_r\bar{R}V_1^r+\partial_z\bar{R}V_1^z)|\tilde{V}_1|^{p-2}r^{p-3}-\frac{N}r\bar{R}_1 N^{\frac{p}2}U_1^\theta|U_1^\theta|^{p-2}r^{p-3}dxdt'|\\
&+\sum_{k\geq2}(kN)^{2\beta_p}|\int_0^t\int_{\mathbb{R}^3}(\partial_r\bar{Q}_kU_k^r+\partial_z\bar{Q}_kU_k^z)|\tilde{U}_k|^{p-2}r^{p-3}+\frac{kN}r\bar{Q}_k (kN)^{\frac{p}2}V_k^\theta|V_k^\theta|^{p-2}r^{p-3}dxdt'|\\
&+\sum_{k\geq2}(kN)^{2\beta_p}|\int_0^t\int_{\mathbb{R}^3}(\partial_r\bar{R}_kV_k^r+\partial_z\bar{R}_kV_k^z)|\tilde{V}_k|^{p-2}r^{p-3}-\frac{kN}r\bar{R}_k (kN)^{\frac{p}2}U_k^\theta|U_k^\theta|^{p-2}r^{p-3}dxdt'|\\
&\leq \frac18(E_p(t)+\mathcal{E}_p(t))+CN^{-\frac{p}{2(p-5)}}D^{\frac{p}2}(t). \end{align*}

\textbf{Lemma 4.3} Let $0<\alpha_p<\frac{p-3-2\beta_p}4$ and $1<\beta_p<\frac{p-1}4$ with $p\in(5,6)$, then for all $t\in[0,T]$
\begin{align*}
N^{-2\alpha_p}|\int_0^t\int_{\mathbb{R}^3}(\phi_1\xi_1|\xi_1|^{p-2}+\psi_1\zeta_1|\zeta_1|^{p-2})r^{p-3}dxdt'|&+\sum_{k\geq2}(kN)^{2\beta_p}|\int_0^t\int_{\mathbb{R}^3}(\phi_k\xi_k|\xi_k|^{p-2}+\psi_k\zeta_k|\zeta_k|^{p-2})r^{p-3}dxdt'|\\
&\leq C(E_p(t)+\mathcal{E}_p(t))^{1+\frac1p}.
\end{align*}

First of all, we admit the above Lemma 4.2 and 4.3. At the end of the section, their proofs will be given. Applying these lemmas, the Proposition 3.1 is obtained.

\textbf{Proof of Proposition 3.1:} Compared with the Navier-Stokes equations, we only deal with the terms relating the temperature function $\eta$. The proof is made up of two parts.

\textbf{Part one: the weight norm of the velocity vector}

According to the work \cite{21}, we know that for any $t\in[0,T]$,
\begin{eqnarray*}
E_p(t)\leq E_p(0)+C(\|U\|_{L_t^\infty(L^3)}+N^{-\frac12})E_p(t)+CE_p^{1+\frac1p}(t)+I_1(t)+I_2(t), \end{eqnarray*}
where
\begin{align*}
I_1(t)=&N^{-2\alpha_p}|\int_0^t\int_{\mathbb{R}^3}\frac1{1+r^2}(\xi_1U_1^z|\tilde{U}_1|^{p-2}+\zeta_1V_1^z|\tilde{V}_1|^{p-2})r^{p-3}dxdt'|\\
&+\sum_{k\geq2}(kN)^{2\beta_p}|\int_0^t\int_{\mathbb{R}^3}\frac1{1+r^2}(\xi_kU_k^z|\tilde{U}_k|^{p-2}+\zeta_kV_k^z|\tilde{V}_k|^{p-2})r^{p-3}dxdt'|,
\end{align*}
and
\begin{align*}
I_2(t)=&N^{-2\alpha_p}|\int_0^t\int_{\mathbb{R}^3}(\partial_rQ_1U_1^r+\partial_zQ_1U_1^z)|\tilde{U}_1|^{p-2}r^{p-3}+\frac{N}rQ_1 N^{\frac{p}2}V_1^\theta|V_1^\theta|^{p-2}r^{p-3}dxdt'|\\
&+N^{-2\alpha_p}|\int_0^t\int_{\mathbb{R}^3}(\partial_rRV_1^r+\partial_zRV_1^z)|\tilde{V}_1|^{p-2}r^{p-3}-\frac{N}rR_1 N^{\frac{p}2}U_1^\theta|U_1^\theta|^{p-2}r^{p-3}dxdt'|\\
&+\sum_{k\geq2}(kN)^{2\beta_p}|\int_0^t\int_{\mathbb{R}^3}(\partial_rQ_kU_k^r+\partial_zQ_kU_k^z)|\tilde{U}_k|^{p-2}r^{p-3}+\frac{kN}rQ_k (kN)^{\frac{p}2}V_k^\theta|V_k^\theta|^{p-2}r^{p-3}dxdt'|\\
&+\sum_{k\geq2}(kN)^{2\beta_p}|\int_0^t\int_{\mathbb{R}^3}(\partial_rR_kV_k^r+\partial_zR_kV_k^z)|\tilde{V}_k|^{p-2}r^{p-3}-\frac{kN}rR_k (kN)^{\frac{p}2}U_k^\theta|U_k^\theta|^{p-2}r^{p-3}dxdt'|.
\end{align*}
In what follows, we shall handle the estimate of $I_j(t),j=1,2$ term by term.

For $I_1(t)$, H\"{o}lder inequality gives
\begin{align*}
\int_{\mathbb{R}^3}|\xi_k| |U_k^z||\tilde{U}_k|^{p-2}r^{p-5}dx&\leq\|r^{\frac{p-5}2}|\tilde{U}_k|^{\frac{p}2}\|_{L^2}^{2-\frac2p}\|r^{\frac{p-5}2}|\xi_k|^{\frac{p}2}\|_{L^2}^{\frac2p}\\
&\lesssim(kN)^{-2}((kN\|r^{\frac{p-5}2}|\tilde{U}_k|^{\frac{p}2}\|_{L^2})^2+(kN\|r^{\frac{p-5}2}|\xi_k|^{\frac{p}2}\|_{L^2})^2).
\end{align*}
Owing to the definition of $E_p(t)$ and $\mathcal{E}_p(t)$, we have
\begin{eqnarray*}
I_1(t)\leq CN^{-2}(E_p(t)+\mathcal{E}_p(t)).
\end{eqnarray*}

As for $I_2(t)$, it is more difficult than $I_1(t)$. According to the incompressible condition, we have the elliptic equations:
\begin{align*}
\mathcal{L}_{kN}Q_k&=-(\partial_r^2+\frac1r\partial_r+\partial_z^2)Q_k+\frac{(kN)^2}{r^2}Q_k\\
&=(\partial_r+\frac1r)(\tilde{U}\cdot\tilde{\nabla}U_k^r+\tilde{U}_k\cdot\tilde{\nabla}U^r-\frac{kN}rU^\theta V_k^r-\frac2rU_k^\theta U^\theta)\\
&-\frac{kN}r(\tilde{U}\cdot\tilde{\nabla}V_k^\theta+\tilde{V}_k\cdot\tilde{\nabla}U^\theta+\frac{kN}rU^\theta U_k^\theta+\frac1r(U^rV_k^\theta+V_k^rU^\theta))\\
&+\partial_z(\tilde{U}\cdot\tilde{\nabla}U_k^z+\tilde{U}_k\cdot\tilde{\nabla}U^z-\frac{kN}rU^\theta V_k^z)+(\partial_r+\frac1r)F_k^r-\frac{kN}rG_k^\theta+\partial_zF_k^z\\
&-(\partial_r+\frac1r)\tilde{\Delta}U_k^r+\frac{kN}r\tilde{\Delta}V_k^\theta-\partial_z\tilde{\Delta}U_k^z-\frac2{r^3}U_k^r-2kN(\partial_r+\frac1r)(\frac1{r^2}V_k^\theta)-\frac1{r^2}\partial_zU_k^z-\frac1{1+r^2}\partial_z\xi_k.
\end{align*}
Write
\begin{eqnarray*}
Q_k(t,r,z)=\tilde{Q}_k(t,r,z)+\bar{Q}_k(t,r,z),
\end{eqnarray*}
where $\bar{Q}_k$ meets the equation:
\begin{eqnarray*}
\mathcal{L}_{kN}\bar{Q}_k=-\frac1{1+r^2}\partial_z\xi_k.
\end{eqnarray*}
Then the rest term $\tilde{Q}_k$ satisfies the Lemma 4.1. Hence, Lemma 4.1 and Lemma 4.2 claims that
\begin{eqnarray*}
I_2(t)\leq (\frac14+C\|U\|_{L_t^\infty(L^3)})(E_p(t)+\mathcal{E}_p(t))+CE_p^{1+\frac1p}(t)+CN^{-2}D^{\frac{p+1}2}
+CN^{-p}(\|U\|_{L_t^\infty(L^3)}+N^{-\frac{p}{2(p-5)}})D^{\frac{p}2}.
\end{eqnarray*}

Summing the estimations, we find that
\begin{align}\label{eq13}
E_p(t)\leq& E_p(0)+C(\|U\|_{L_t^\infty(L^3)}+N^{-\frac12})(E_p(t)+\mathcal{E}_p(t))+CE_p^{1+\frac1p}(t)+CN^{-2}D^{\frac{p+1}2}(t)\notag\\
&+CN^{-p}(\|U\|_{L_t^\infty(L^3)}+N^{-\frac{p}{2(p-5)}})D^{\frac{p}2}(t). \end{align}

\textbf{Part two: the weight norm of the temperature functions}

Firstly, we estimate the weight norm of the temperature functions $\zeta_k$. For $p\in(5,6)$, by taking $L^2$ inner product of $\zeta_k$ equations \eqref{eq8} with $r^{p-3}|\zeta_k|^{p-2}\zeta_k$, we obtain
\begin{align*}
\int_{\mathbb{R}^3}&\left(\frac12\frac{d}{dt}(\zeta_k)^2-\tilde{\Delta}\zeta_k\cdot \zeta_k\right)|\zeta_k|^{p-2}r^{p-3}dx+(kN)^2\|r^{\frac{p-5}2}|\zeta_k|^{\frac{p}2}\|_{L^2}^2+\int_{\mathbb{R}^3}\tilde{L}\zeta_k\cdot\zeta_k|\zeta_k|^{p-2}r^{p-3}dx\\
&=-\int_{\mathbb{R}^3}\left(\tilde{U}\cdot\tilde{\nabla}\zeta_k+\tilde{V}_k\cdot\tilde{\nabla}\xi+\frac{kN}rU^\theta \xi_k-\frac{2r}{1+r^2}U^r\zeta_k-\frac{2r}{1+r^2}V_k^r\xi+\psi_k\right)\zeta_k|\zeta_k|^{p-2}r^{p-3}dx.
\end{align*}
Using integration by parts, we have
\begin{align*}
-\int_{\mathbb{R}^3}\tilde{\Delta}\zeta_k\cdot \zeta_k|\zeta_k|^{p-2}r^{p-3}dx&=-2\pi\int_{\mathbb{R}^+\times\mathbb{R}}(\partial_r^2+\frac1r\partial_r+\partial_z^2)\zeta_k\cdot \zeta_k|\zeta_k|^{p-2}r^{p-2}drdz\\
&=2\pi\int_{\mathbb{R}^+\times\mathbb{R}}\left(|\tilde{\nabla}\zeta_k|^2+\frac{p-2}4\tilde{\nabla}(\zeta_k)^2\cdot |\zeta_k|^{-2}\tilde{\nabla}|\zeta_k|^2+\frac{p-3}{2r}\partial_r(\zeta_k)^2\right)|\zeta_k|^{p-2}r^{p-2}drdz\\
&=\int_{\mathbb{R}^3}|\tilde{\nabla}\zeta_k|^2|\tilde{V}_k|^{p-2}r^{p-3}dx+(p-2)\int_{\mathbb{R}^3}|\tilde{\nabla}|\zeta_k||^2|\zeta_k|^{p-2}r^{p-3}dx\\
&~~~~~~~~+\pi(p-3)\int_0^\infty\int_{\mathbb{R}}\partial_r|\zeta_k|^2\cdot|\zeta_k|^{p-2}r^{p-3}drdz\\ &\geq\int_{\mathbb{R}^3}|\tilde{\nabla}\zeta_k|^2|\tilde{V}_k|^{p-2}r^{p-3}dx-2\int_{\mathbb{R}^3}|\zeta_k|^pr^{p-5}dx,
\end{align*}
where in the last step, we used integration by part and $p\in(5,6)$ so that
\begin{align*}
\pi(p-3)\int_0^\infty\int_{\mathbb{R}}\partial_r|\zeta_k|^2\cdot|\zeta_k|^{p-2}r^{p-3}drdz
&=-\frac{2\pi(p-3)^2}p\int_0^\infty\int_{\mathbb{R}}|\zeta_k|^pr^{p-4}drdz+\frac{2\pi(p-3)}p\int_{\mathbb{R}}\lim_{r\rightarrow\infty}(r^{p-3}|\zeta_k|^p)dz\\
&\geq-2\int_{\mathbb{R}^3}|\zeta_k|^pr^{p-5}dx.
\end{align*}
In addition, for the linear term $\tilde{L}\zeta_k$, against by using integration by parts, we have
\begin{align*}
\int_{\mathbb{R}^3}(\tilde{L}\zeta_k)\zeta_k|\zeta_k|^{p-2}r^{p-3}dx&=2\pi\int_{\mathbb{R}^+\times\mathbb{R}}(\frac{4r}{1+r^2}\partial_r+\frac{4(1-r^2)}{(1+r^2)^2})\zeta_k\cdot \zeta_k|\zeta_k|^{p-2}r^{p-2}drdz\\
&\geq-8\pi\int_0^\infty\int_{\mathbb{R}}|\zeta_k|^pr^{p-4}drdz-\frac{8\pi}p\int_0^\infty\int_{\mathbb{R}}\frac{r^2(p-1+(p-3)r^2)}{(1+r^2)^2}|\zeta_k|^pr^{p-4}drdz\\
&~~~~~+\frac{8\pi}p\int_{\mathbb{R}}\lim_{r\rightarrow\infty}(\frac{r^{p-1}}{1+r^2}|\zeta_k|^p)dz\\ &\geq-8\int_{\mathbb{R}^3}|\zeta_k|^pr^{p-5}dx,
\end{align*}
where in the last step, we used the inequality for all $r\geq0$ and $p\in(5,6)$
\begin{eqnarray*}
|\frac{(1-r^2)r^2}{(1+r^2)^2}|\leq1,~~~\frac{r^2(p-1+(p-3)r^2)}{(1+r^2)^2}\leq p-1.
\end{eqnarray*}
Then we get, that
\begin{eqnarray*}
\frac1p\frac{d}{dt}\|r^{\frac{p-3}2}|\zeta_k|^{\frac{p}2}\|_{L^2}^2+\|r^{\frac{p-3}2}|\tilde{\nabla}\zeta_k||\zeta_k|^{\frac{p}2-1}\|_{L^2}^2
+((kN)^2-10)\|r^{\frac{p-5}2}|\zeta_k|^{\frac{p}2}\|_{L^2}^2 \leq\sum_{j=1}^4|II_k^j(t)|,
\end{eqnarray*}
where
\begin{align*}
II_k^1(t)&=\int_{\mathbb{R}^3}(\tilde{U}\cdot\tilde{\nabla}\zeta_k)\cdot \zeta_k|\zeta_k|^{p-2}r^{p-3}dx,~~II_k^2(t)=\int_{\mathbb{R}^3}(\tilde{V}_k\cdot\tilde{\nabla}\xi)\cdot \zeta_k|\zeta_k|^{p-2}r^{p-3}dx,\\
II_k^3(t)&=\int_{\mathbb{R}^3}(kNU^\theta \xi_k-\frac{2r^2}{1+r^2}U^r\zeta_k-\frac{2r^2}{1+r^2}V_k^r\xi)\zeta_k|\zeta_k|^{p-2}r^{p-4}dx,~~~
II_k^4(t)=\int_{\mathbb{R}^3}\psi_k\zeta_k|\zeta_k|^{p-2}r^{p-3}dx.
\end{align*}
In what follows, we shall handle the estimate of $II_k^j(t)$ term by term.

Firstly, since $r^{\frac{p-3}2}|\zeta_k|^{\frac{p}2}$ is axisymmetric, it follows from Sobolev embedding theorem that
\begin{eqnarray*}
\|r^{\frac{p-3}2}|\zeta_k|^{\frac{p}2}\|_{L^6}\lesssim\|\tilde{\nabla}(r^{\frac{p-3}2}|\zeta_k|^{\frac{p}2})\|_{L^2}\lesssim\|r^{\frac{p-3}2}\tilde{\nabla}\zeta_k|\zeta_k|^{\frac{p}2-1}\|_{L^2}+\|r^{\frac{p-5}2}|\zeta_k|^{\frac{p}2}\|_{L^2},
\end{eqnarray*}
from which, we infer that
\begin{eqnarray*}
|II_k^1(t)|\lesssim\|\tilde{U}\|_{L^3}\|r^{\frac{p-3}2}|\tilde{\nabla}\zeta_k||\zeta_k|^{\frac{p}2-1}\|_{L^2}\|r^{\frac{p-3}2}|\zeta_k|^{\frac{p}2}\|_{L^6}
\lesssim\|\tilde{U}\|_{L^3}(\|r^{\frac{p-3}2}|\tilde{\nabla}\zeta_k||\zeta_k|^{\frac{p}2-1}\|_{L^2}^2
+\|r^{\frac{p-5}2}|\zeta_k|^{\frac{p}2}\|_{L^2}^2).
\end{eqnarray*}
For $II_k^2(t)$, we get, by using integration by parts and the incompressible condition, that
\begin{align*}
II_k^2(t)&=2\pi\int_0^\infty\int_{\mathbb{R}}(\tilde{V}_k\cdot\tilde{\nabla}\xi)\cdot \zeta_k|\zeta_k|^{p-2}r^{p-2}drdz\\
&=-\int_{\mathbb{R}^3}(\xi\cdot(\tilde{V}_k\cdot\tilde{\nabla}\zeta_k)+\xi\zeta_k)((p-3)\frac{V_k^r}r-\frac{kN}rU_k^\theta))|\zeta_k|^{p-2}r^{p-3}dx-(p-2)\int_{\mathbb{R}^3}\xi\zeta_k(\tilde{V}_k\cdot\nabla|\zeta_k|)|\zeta_k|^{p-3}r^{p-3}dx.
\end{align*}
By applying H\"{o}lder inequality and Young inequality, we find
\begin{align*}
|II_k^2(t)|&\lesssim\|\xi\|_{L^3}\|r^{\frac{p-3}2}|\zeta_k|^{\frac{p}2}\|_{L^6}^{1-\frac2p}\|r^{\frac{p-3}2}|\tilde{V}_k|^{\frac{p}2}\|_{L^6}^{\frac2p}\|r^{\frac{p-3}2}|\tilde{\nabla}\zeta_k||\zeta_k|^{\frac{p}2-1}\|_{L^2}+\|\xi\|_{L^3}\|r^{\frac{p-3}2}|\zeta_k|^{\frac{p}2}\|_{L^6}\|r^{\frac{p-5}2}|\zeta_k|^{\frac{p}2}\|_{L^2}^{1-\frac2p}\|r^{\frac{p-5}2}|\tilde{V}_k|^{\frac{p}2}\|_{L^2}^{\frac2p}\\
&+\|\xi\|_{L^3}\|r^{\frac{p-3}2}|\zeta_k|^{\frac{p}2}\|_{L^6}(kN\|r^{\frac{p-5}2}|\zeta_k|^{\frac{p}2}\|_{L^2})^{1-\frac2p}(kN\|r^{\frac{p-5}2}|U_k^\theta|^{\frac{p}2}\|_{L^2})^{\frac2p}\\
&\lesssim\|\xi\|_{L^3}(\|r^{\frac{p-3}2}|\tilde{\nabla}\tilde{V}_k||\tilde{V}_k|^{\frac{p}2-1}\|_{L^2}^2+\|r^{\frac{p-3}2}|\tilde{\nabla}\zeta_k||\zeta_k|^{\frac{p}2-1}\|_{L^2}^2\\
&+(kN)^2\|r^{\frac{p-5}2}|\tilde{V}_k|^{\frac{p}2}\|_{L^2}^2+(kN)^2\|r^{\frac{p-5}2}|\zeta_k|^{\frac{p}2}\|_{L^2}^2+(kN)^{\frac{p}2+2}\|r^{\frac{p-5}2}|U_k^\theta|^{\frac{p}2}\|_{L^2}^2).
\end{align*}
At same time, due to $N\geq1$, $II_k^3(t)$ also can be handled as follows
\begin{align*}
|II_k^3(t)|&\leq \|U^\theta\|_{L^3}\|r^{\frac{p-3}2}|\zeta_k|^{\frac{p}2}\|_{L^6}(kN\|r^{\frac{p-5}2}|\zeta_k|^{\frac{p}2}\|_{L^2})^{1-\frac2p}(kN\|r^{\frac{p-5}2}|\xi_k|^{\frac{p}2}\|_{L^2})^{\frac2p}\\
&+~\|U^r\|_{L^3}\|r^{\frac{p-3}2}|\zeta_k|^{\frac{p}2}\|_{L^6}\|r^{\frac{p-5}2}|\zeta_k|^{\frac{p}2}\|_{L^2}+\|\xi\|_{L^3}\|r^{\frac{p-3}2}|\zeta_k|^{\frac{p}2}\|_{L^6}\|r^{\frac{p-5}2}|\zeta_k|^{\frac{p}2}\|_{L^2}^{1-\frac2p}\|r^{\frac{p-5}2}|V_k^r|^{\frac{p}2}\|_{L^2}^{\frac2p}\\
&\lesssim(\|U\|_{L^3}+\|\xi\|_{L^3})(\|r^{\frac{p-3}2}|\tilde{\nabla}\zeta_k||\zeta_k|^{\frac{p}2-1}\|_{L^2}^2+(kN)^2\|r^{\frac{p-5}2}|\zeta_k|^{\frac{p}2}\|_{L^2}^2\\
&~~+(kN)^2\|r^{\frac{p-5}2}|\tilde{V}_k|^{\frac{p}2}\|_{L^2}^2+(kN)^2\|r^{\frac{p-5}2}|\xi_k|^{\frac{p}2}\|_{L^2}^2).
\end{align*}
By substituting these estimates and multiplying the inequalities by
$N^{-2\alpha_p}$ for $k=1$ and by $(kN)^{2\beta_p}$ for $k\geq2$, and then integrating them with respect to time over $[0,t]$, and finally summing up the resulting inequalities for $k\in\mathbb{N}^+$, we achieve by Lemma 4.3
\begin{align}\label{eq14}
&\frac1pN^{-2\alpha_p}\|r^{\frac{p-3}2}|\zeta_1(t)|^{\frac{p}2}\|_{L^2}^2+\frac1p\sum_{k=2}^\infty(kN)^{2\beta_p}\|r^{\frac{p-3}2}|\zeta_k(t)|^{\frac{p}2}\|_{L^2}^2+N^{-2\alpha_p}\int_0^t\|r^{\frac{p-3}2}|\tilde{\nabla}\zeta_1||\zeta_1|^{\frac{p}2-1}\|_{L^2}^2dt'\notag\\
&+\sum_{k=2}^\infty(kN)^{2\beta_p}\int_0^t\|r^{\frac{p-3}2}|\tilde{\nabla}\zeta_k||\zeta_k|^{\frac{p}2-1}\|_{L^2}^2dt'+N^{-2\alpha_p}(N^2-10)\int_0^t\|r^{\frac{p-5}2}|\zeta_1|^{\frac{p}2}\|_{L^2}^2dt'\notag\\
&+\sum_{k=2}^\infty(kN)^{2\beta_p}((kN)^2-10)\int_0^t\|r^{\frac{p-5}2}|\zeta_k|^{\frac{p}2}\|_{L^2}^2dt'\leq\frac1p(N^{-2\alpha_p}\|r^{\frac{p-3}2}|\zeta_1(0)|^{\frac{p}2}\|_{L^2}^2+\sum_{k=2}^\infty(kN)^{2\beta_p}\|r^{\frac{p-3}2}|\zeta_k(0)|^{\frac{p}2}\|_{L^2}^2)\notag\\
&+C(\|U\|_{L_t^\infty(L^3)}+\|\xi\|_{L_t^\infty(L^3)})\mathcal{E}_p(t)+C(E_p(t)+\mathcal{E}_p(t))^{1+\frac1p}.
\end{align}

Meanwhile, for the another temperature function $\xi_k$, we also obtain that
\begin{align}\label{eq15}
&\frac1pN^{-2\alpha_p}\|r^{\frac{p-3}2}|\xi_1(t)|^{\frac{p}2}\|_{L^2}^2+\frac1p\sum_{k=2}^\infty(kN)^{2\beta_p}\|r^{\frac{p-3}2}|\xi_k(t)|^{\frac{p}2}\|_{L^2}^2+N^{-2\alpha_p}\int_0^t\|r^{\frac{p-3}2}|\tilde{\nabla}\xi_1||\xi_1|^{\frac{p}2-1}\|_{L^2}^2dt'\notag\\
&+\sum_{k=2}^\infty(kN)^{2\beta_p}\int_0^t\|r^{\frac{p-3}2}|\tilde{\nabla}\xi_k||\xi_k|^{\frac{p}2-1}\|_{L^2}^2dt'+N^{-2\alpha_p}(N^2-10)\int_0^t\|r^{\frac{p-5}2}|\xi_1|^{\frac{p}2}\|_{L^2}^2dt'\notag\\
&+\sum_{k=2}^\infty(kN)^{2\beta_p}((kN)^2-10)\int_0^t\|r^{\frac{p-5}2}|\xi_k|^{\frac{p}2}\|_{L^2}^2dt'
\leq\frac1p(N^{-2\alpha_p}\|r^{\frac{p-3}2}|\xi_1(0)|^{\frac{p}2}\|_{L^2}^2+\sum_{k=2}^\infty(kN)^{2\beta_p}\|r^{\frac{p-3}2}|\xi_k(0)|^{\frac{p}2}\|_{L^2}^2)\notag\\
&+C(\|U\|_{L_t^\infty(L^3)}+\|\xi\|_{L_t^\infty(L^3)})\mathcal{E}_p(t)+C(E_p(t)+\mathcal{E}_p(t))^{1+\frac1p}.
\end{align}

At last, combining \eqref{eq13}-\eqref{eq15}, we have that for all $0\leq t\leq T$,
\begin{align*}
E_p(t)+\mathcal{E}_p(t)&\leq E_p(0)+\mathcal{E}_p(0)+(\frac14+C(\|U\|_{L_t^\infty(L^3)}+\|\xi\|_{L_t^\infty(L^3)}+\frac1{\sqrt{N}}))(E_p(t)+\mathcal{E}_p(t))\\
&+C(E_p(t)+\mathcal{E}_p(t))^{1+\frac1p}+CN^{-2}D^{\frac{p+1}2}(t)
+CN^{-p}(\|U\|_{L_t^\infty(L^3)}+N^{-\frac{p}{2(p-5)}})D^{\frac{p}2}(t). \end{align*}
Choosing $N$ such that $C/\sqrt{N}\leq1/4$, then we complete the proof. $\square$

In the end, we will prove the Lemma 4.2 and Lemma 4.3. The Lemma 4.2 is a application of Proposition 4.1. But the Lemma 4.3 is more complicated.

\textbf{The proof of Lemma 4.2}: Here we shall only present the estimate of $\bar{R}_k$ part. The rest parts can be handled exactly along the same line.
Due to the assumption in Lemma 4.2, we know that
\begin{eqnarray*}
\bar{R}_k=\mathcal{L}_{kN}^{-1}\Psi_k,~~~\Psi_k=-\frac1{1+r^2}\zeta_k.
\end{eqnarray*}
The proposition 4.1 shows for any small enough $\epsilon>0$,
\begin{align*}
\|r^{\frac{p-3}2}\mathcal{L}_{kN}^{-1}\Psi_k\|_{L_t^{\frac{10}{7-p}}(L^{\frac{10}{7-p}})}&\lesssim(kN)^{-(\frac65-\frac{p}{10}-\frac1p)}\|\frac{r^2}{1+r^2}\zeta_kr^{\frac{p-2}5-\frac3p}\|_{L_t^{\frac{10}{7-p}}(L^p)}\\
&\lesssim(kN)^{-(\frac65-\frac{p}{10}-\frac1p-\epsilon)}\|r^{\frac{p-3}2}|\zeta_k|^{\frac{p}2}\|_{L_t^\infty(L^2)}^{\frac2p-\frac{7-p}5}\|r^{\frac{p-5}2}|\zeta_k|^{\frac{p}2}\|_{L_t^2(L^2)}^{\frac{7-p}5}.
\end{align*}
from which, we infer
\begin{align*}
&\left|\int_0^t\int_{\mathbb{R}^3}\tilde{\nabla}\mathcal{L}_{kN}^{-1}\Psi_k\cdot|\tilde{U}_k|^{p-2}\tilde{V}_kr^{p-3}dxdt'\right|\lesssim\|r^{\frac{p-3}2}\mathcal{L}_{kN}^{-1}\Psi_k\|_{L_t^{\frac{10}{7-p}}(L^{\frac{10}{7-p}})}\\
&~~~~~~\times(\|r^{\frac{p-3}2}|\tilde{\nabla}\tilde{V}_k||\tilde{V}_k|^{\frac{p}2-1}\|_{L_t^2(L^2)}+\|r^{\frac{p-5}2}|\tilde{V}_k|^{\frac{p}2}\|_{L_t^2(L^2)})
\|r^{\frac{p-5}2}|\tilde{V}_k|^{\frac{p}2}\|_{L_t^2(L^2)}^{\frac35}\|r^{-\frac32}\tilde{V}_k\|_{L_t^2(L^2)}^{\frac{p}5-1}\\
&~~~~~~\lesssim(kN)^{-(\frac{16}5-\frac{3p}{10}-\frac1p-\epsilon)}\|r^{\frac{p-3}2}|\zeta_k|^{\frac{p}2}\|_{L_t^\infty(L^2)}^{\frac2p-\frac{7-p}5}(kN\|r^{\frac{p-5}2}|\zeta_k|^{\frac{p}2}\|_{L_t^2(L^2)})^{\frac{7-p}5}
\|r^{\frac{p-3}2}|\tilde{V}_k|^{\frac{p}2}\|_{L_t^\infty(L^2)}\\
&~~~~~~\times(\|r^{\frac{p-3}2}|\tilde{\nabla}\tilde{V}_k||\tilde{V}_k|^{\frac{p}2-1}\|_{L_t^2(L^2)}+\|r^{\frac{p-5}2}|\tilde{V}_k|^{\frac{p}2}\|_{L_t^2(L^2)})(kN\|r^{\frac{p-5}2}|\tilde{V}_k|^{\frac{p}2}\|_{L_t^2(L^2)})^{\frac35}\|r^{-\frac32}\tilde{V}_k\|_{L_t^2(L^2)}^{\frac{p}5-1}.
\end{align*}
Observing that for all $p\in(5,6)$
\begin{eqnarray*}
(\frac2p-\frac{7-p}5)+1+(2-\frac{p}5)+(\frac{p}5-1)>2,
\end{eqnarray*}
then we get
\begin{align*}
\left|\int_0^t\int_{\mathbb{R}^3}\tilde{\nabla}\mathcal{L}_{kN}^{-1}\Psi_k\cdot|\tilde{V}_k|^{p-2}\tilde{V}_kr^{p-3}dxdt'\right|
\lesssim\sup_{k\in\mathbb{N}^+}(kN)^{-(\frac{16}5-\frac{3p}{10}-\frac1p-\epsilon)}\cdot (E_p(t)+\mathcal{E}_p(t))^{\frac45+\frac1p}(\sum_{k=1}^\infty\|r^{-\frac32}\tilde{V}_k\|_{L_t^2(L^2)}^2)^{\frac{p}{10}-\frac12}.
\end{align*}
The requirement of $\beta_p$, ones have
\begin{eqnarray*}
\frac1p+\frac{3p}{10}-\frac{16}5+2\beta_p(\frac15-\frac1p)<\frac1p+\frac{3p}{10}-\frac{16}5+\frac{p-1}2(\frac15-\frac1p)<-1.15.
\end{eqnarray*}
Hence, taking $\epsilon=0.05$, Young inequality claims that
\begin{align*}
\left|\int_0^t\int_{\mathbb{R}^3}\tilde{\nabla}\mathcal{L}_{kN}^{-1}\Psi_k\cdot|\tilde{V}_k|^{p-2}\tilde{V}_kr^{p-3}dxdt'\right|&\leq CN^{-\frac1{10}}(E_p(t)+\mathcal{E}_p(t))^{\frac45+\frac1p}(\sum_{k=1}^\infty\|r^{-\frac32}\tilde{V}_k\|_{L_t^2(L^2)}^2)^{\frac{p}{10}-\frac12}\\
&\leq\frac18(E_p(t)+\mathcal{E}_p(t))+CN^{-\frac{p}{2(p-5)}}(\sum_{k=1}^\infty\|r^{-\frac32}\tilde{V}_k\|_{L_t^2(L^2)}^2)^{\frac{p}2}.
\end{align*}
Summing theses estimations, we compete the proof. $\square$

\textbf{The proof of Lemma 4.3:} Here we shall only present the estimate of $\psi_k\zeta_k$ part. The rest parts can be handled exactly along the same line. Decomposed $\psi_k$ as
\begin{eqnarray*}
\psi_k=\psi_k^1+\psi_k^2,
\end{eqnarray*}
where
\begin{align*}
\psi_k^1&=\frac12\sum_{k_1+k_2=k}\left(\frac{k_2N}rU_{k_1}^\theta \zeta_{k_2}+\tilde{V}_{k_1}\cdot\tilde{\nabla}\zeta_{k_2}-\frac{2r}{1+r^2}V_{k_1}^r\zeta_{k_2}\right)\\
&~~~~~~~~~~~~~~~~~~~~~+\frac12\sum_{|k_1-k_2|=k}\left(-\frac{k_2N}rU_{k_1}^\theta \zeta_{k_2}+\tilde{V}_{k_1}\cdot\tilde{\nabla}\zeta_{k_2}-\frac{2r}{1+r^2}V_{k_1}^r\zeta_{k_2}\right),\\
\psi_k^2&=\frac12\sum_{k_1+k_2=k}\left(-\tilde{U}_{k_1}\cdot\tilde{\nabla}\xi_{k_2}+\frac{k_2N}rV_{k_1}^\theta \xi_{k_2}+\frac{2r}{1+r^2}U_{k_1}^r\xi_{k_2}\right)\\
&~~~~~~~~~~~~~~~~~~~~~+\frac12\sum_{|k_1-k_2|=k}\left(\tilde{U}_{k_1}\cdot\tilde{\nabla}\xi_{k_2}+\frac{k_2N}rV_{k_1}^\theta \xi_{k_2}-\frac{2r}{1+r^2}U_{k_1}^r\xi_{k_2}\right).
\end{align*}
Because of that the structure of $\psi_k^1$ is same as $\psi_k^2$, we only deal with $\psi_k^1$.

In all that follows, we always denote
\begin{eqnarray*}
\Omega_k=\{(k_1,k_2)\in\mathbb{N}^+\times\mathbb{N}^+: |k_1-k_2|=k~or~k_1+k_2=k\},~~~\mathcal{H}_p(t)=E_p(t)+\mathcal{E}_p(t).
\end{eqnarray*}
By using integration by parts, ones get
\begin{align*}
&\sum_{(k_1,k_2)\in\Omega_k}|\int_{\mathbb{R}^3}(\tilde{V}_{k_1}\cdot\tilde{\nabla}\zeta_{k_2})\zeta_k|\zeta_k|^{p-2}r^{p-3}dx|\\
&=\sum_{(k_1,k_2)\in\Omega_k}|\int_{\mathbb{R}^3}\zeta_{k_2}(\tilde{V}_{k_1}\cdot\tilde{\nabla}\zeta_k+(p-2)\zeta_k|\zeta_k|^{-1}\tilde{V}_{k_1}\cdot\tilde{\nabla}|\zeta_k|+\frac{\zeta_k}r((p-3)V_{k_1}^r-k_1NU_{k_1}^\theta))|\zeta_k|^{p-2}r^{p-3}dx|\\
&\lesssim\sum_{(k_1,k_2)\in\Omega_k}\int_{\mathbb{R}^3}|\zeta_{k_2}|(|\tilde{V}_{k_1}|(|\tilde{\nabla}\zeta_k|+\frac1r|\zeta_k|)+\frac{k_1N}r|U_{k_1}^\theta||\zeta_k|)|\zeta_k|^{p-2}r^{p-3}dx|,
\end{align*}
with the incompressible conditions. So, combing with $\frac{2r}{1+r^2}\leq\frac2r$, we decompose it as
\begin{eqnarray}\label{eq16}
|\int_{\mathbb{R}^3}\psi_k^1\zeta_k|\zeta_k|^{p-2}r^{p-3}dx|\lesssim\sum_{j=1}^2III_k^j(t),
\end{eqnarray}
where
\begin{align*}
III_k^1(t)=\sum_{(k_1,k_2)\in\Omega_k}\int_{\mathbb{R}^3}|\tilde{V}_{k_1}||\zeta_{k_2}|(|\tilde{\nabla}\zeta_k|+\frac1r|\zeta_k|)|\zeta_k|^{p-2}r^{p-3}dx,~~III_k^2(t)=\sum_{(k_1,k_2)\in\Omega_k}\int_{\mathbb{R}^3}(k_1+k_2)N|U_{k_1}^\theta ||\zeta_{k_2}||\tilde{V}_k|^{p-1}r^{p-4}dx.
\end{align*}
Below we shall handle the estimate term by term above.

\textbf{The estimate of $III_k^1$ when $k=1$:}

There is no $k_1,k_2\in\mathbb{N}^+$ satisfying $k_1+k_2=1$. So that in this case, we have
\begin{align*}
III_1^1(t)&=\sum_{j=1}^\infty\int_{\mathbb{R}^3}|\tilde{V}_j||\zeta_{j+1}|(|\tilde{\nabla}\zeta_1|+\frac1r|\zeta_1|)|\zeta_1|^{p-2}r^{p-3}dx+\sum_{j=1}^\infty\int_{\mathbb{R}^3}|\tilde{V}_{j+1}||\zeta_j|(|\tilde{\nabla}\zeta_1|+\frac1r|\zeta_1|)|\zeta_1|^{p-2}r^{p-3}dx\\
&:=III_1^{1,1}(t)+III_1^{1,2}(t).
 \end{align*}
Then for $III_1^{1,1}(t)$, we get, by using H\"{o}lder's inequality, that
\begin{align*}
N^{-2\alpha_p}III_1^{1,1}&\lesssim(N^{1-\alpha_p}\|r^{\frac{p-5}2}|\tilde{V}_1|^{\frac{p}2}\|_{L^2})^{\frac12-\frac3{2p}}
(N^{-\alpha_p}\|r^{\frac{p-3}2}|\tilde{V}_1|^{\frac{p}2}\|_{L^6})^{\frac7{2p}-\frac12}((2N)^{1+\beta_p}\|r^{\frac{p-5}2}|\zeta_2|^{\frac{p}2}\|_{L^2})^{\frac12-\frac3{2p}}\\
&\times((2N)^{\beta_p}\|r^{\frac{p-3}2}|\zeta_2|^{\frac{p}2}\|_{L^6})^{\frac7{2p}-\frac12}(N^{-\alpha_p}\|r^{\frac{p-3}2}|\tilde{\nabla}\zeta_1||\zeta_1|^{\frac{p}2-1}\|_{L^2}+N^{1-\alpha_p}\|r^{\frac{p-5}2}|\zeta_1|^{\frac{p}2}\|_{L^2})\\
&\times(N^{-\alpha_p}\|r^{\frac{p-3}2}|\zeta_1|^{\frac{p}2}\|_{L^2})^{\frac2p}
(N^{-\alpha_p}\|r^{\frac{p-3}2}|\zeta_1|^{\frac{p}2}\|_{L^6})^{1-\frac4p}N^{-1+\frac3p-\frac2p\beta_p}\\
&+\sum_{j=2}^\infty((jN)^{1+\beta_p}\|r^{\frac{p-5}2}|\tilde{V}_j|^{\frac{p}2}\|_{L^2})^{\frac12-\frac3{2p}}
((jN)^{\beta_p}\|r^{\frac{p-3}2}|\tilde{V}_j|^{\frac{p}2}\|_{L^6})^{\frac7{2p}-\frac12}\\
&\times(((j+1)N)^{1+\beta_p}\|r^{\frac{p-5}2}|\zeta_{j+1}|^{\frac{p}2}\|_{L^2})^{\frac12-\frac3{2p}}
(((j+1)N)^{\beta_p}\|r^{\frac{p-3}2}|\zeta_{j+1}|^{\frac{p}2}\|_{L^6})^{\frac7{2p}-\frac12}\\
&\times(N^{-\alpha_p}\|r^{\frac{p-3}2}|\zeta_1|^{\frac{p}2}\|_{L^2}+N^{1-\alpha_p}\|r^{\frac{p-5}2}|\zeta_1|^{\frac{p}2}\|_{L^2})\\
&\times(N^{-\alpha_p}\|r^{\frac{p-3}2}|\zeta_1|^{\frac{p}2}\|_{L^2})^{\frac2p}
(N^{\frac{p}4-\alpha_p}\|r^{\frac{p-3}2}|\zeta_1|^{\frac{p}2}\|_{L^6})^{1-\frac4p}N^{-\frac2p\alpha_p}(jN)^{-1+\frac3p-\frac4p\beta_p},
\end{align*}
where we have used the fact that $\frac{j+1}j\in(1,\frac32]$ for every $j\geq2$. Observing that for $5<p<6$ and $\beta_p>1$, one has
\begin{eqnarray*}
(-1+\frac3p-\frac4p\beta_p)\frac{p}{p-2}<-\frac{p+1}{p-2}<-\frac74
\Rightarrow\sum_{j=2}^\infty j^{(-1+\frac3p-\frac4p\beta_p)\frac{p}{p-2}}<+\infty.
\end{eqnarray*}
Then we get, by integrating in time over $[0,t]$, that
\begin{align*}
N^{-2\alpha_p}\int_0^tIII_1^{1,1}(t')dt'&\lesssim N^{-1+\frac3p-\frac2p\beta_p}\mathcal{H}_p^{1+\frac1p}(t)+N^{-1+\frac3p-\frac4p\beta_p-\frac2p\alpha_p}\tilde{H}_p^{1+\frac1p}(t)\left(\sum_{j=2}^\infty j^{(-1+\frac3p-\frac4p\beta_p)\frac{p}{p-2}}\right)^{1-\frac2p}\\
&\lesssim N^{-1+\frac3p-\frac2p\beta_p}\mathcal{H}_p^{1+\frac1p}(t).
\end{align*}

As for $III_1^{1,2}(t)$, we get, by using H\"{o}lder's inequality, that
\begin{align*}
N^{-2\alpha_p}III_1^{1,2}&\lesssim((2N)^{1+\beta_p}\|r^{\frac{p-5}2}|\tilde{V}_2|^{\frac{p}2}\|_{L^2})^{\frac12-\frac3{2p}}
((2N)^{\beta_p}\|r^{\frac{p-3}2}|\tilde{V}_2|^{\frac{p}2}\|_{L^6})^{\frac7{2p}-\frac12}(N^{1-\alpha_p}\|r^{\frac{p-5}2}|\zeta_1|^{\frac{p}2}\|_{L^2})^{\frac12-\frac3{2p}}\\
&\times(N^{-\alpha_p}\|r^{\frac{p-3}2}|\zeta_1|^{\frac{p}2}\|_{L^6})^{\frac7{2p}-\frac12}(N^{-\alpha_p}\|r^{\frac{p-3}2}|\zeta_1|^{\frac{p}2}\|_{L^2}+N^{1-\alpha_p}\|r^{\frac{p-5}2}|\zeta_1|^{\frac{p}2}\|_{L^2})\\
&\times(N^{-\alpha_p}\|r^{\frac{p-3}2}|\zeta_1|^{\frac{p}2}\|_{L^2})^{\frac2p}
(N^{-\alpha_p}\|r^{\frac{p-3}2}|\zeta_1|^{\frac{p}2}\|_{L^6})^{1-\frac4p}N^{-1+\frac3p-\frac2p\beta_p}\\
&+\sum_{j=2}^\infty(((j+1)N)^{1+\beta_p}\|r^{\frac{p-5}2}|\tilde{V}_{j+1}|^{\frac{p}2}\|_{L^2})^{\frac12-\frac3{2p}}
(((j+1)N)^{\beta_p}\|r^{\frac{p-3}2}|\tilde{V}_{j+1}|^{\frac{p}2}\|_{L^6})^{\frac7{2p}-\frac12}\\
&\times((jN)^{1+\beta_p}\|r^{\frac{p-5}2}|\zeta_j|^{\frac{p}2}\|_{L^2})^{\frac12-\frac3{2p}}
((jN)^{\beta_p}\|r^{\frac{p-3}2}|\zeta_j|^{\frac{p}2}\|_{L^6})^{\frac7{2p}-\frac12}\\
&\times(N^{-\alpha_p}\|r^{\frac{p-3}2}|\zeta_1|^{\frac{p}2}\|_{L^2}+N^{1-\alpha_p}\|r^{\frac{p-5}2}|\zeta_1|^{\frac{p}2}\|_{L^2})\\
&\times(N^{-\alpha_p}\|r^{\frac{p-3}2}|\zeta_1|^{\frac{p}2}\|_{L^2})^{\frac2p}
(N^{-\alpha_p}\|r^{\frac{p-3}2}|\zeta_1a|^{\frac{p}2}\|_{L^6})^{1-\frac4p}N^{-\frac2p\alpha_p}(jN)^{-1+\frac3p-\frac4p\beta_p},
\end{align*}
where we have used the fact that $\frac{j+1}j\in(1,\frac32]$ for every $j\geq2$. The same argument gives
\begin{eqnarray*}
N^{-2\alpha_p}\int_0^tIII_1^{1,2}(t')dt'\lesssim N^{-1+\frac3p-\frac2p\beta_p}\mathcal{H}_p^{1+\frac1p}(t).
\end{eqnarray*}

\textbf{The estimate of $III_k^1$ when $k=2$:}

The estimate of $III_2^1$ is similar to that of $III_1^1$, except now the $k_1+k_2=2$ part in the summation is nontrivial. Indeed it allows $k_1=k_2=1$, which corresponds to the term
\begin{align*}
(2N)^{2\beta_p}\int_0^t\int_{\mathbb{R}^3}|\tilde{V}_1||\zeta_1|(|\tilde{\nabla}\zeta_2|&+\frac1r|\zeta_2|)|\zeta_2|^{p-2}r^{p-3}dxdt'\\
&\lesssim N^{-1+\frac3p+\frac2p\beta_p+\frac4p\alpha_p}\int_0^t(N^{1-\alpha_p}\|r^{\frac{p-5}2}|\tilde{V}_1|^{\frac{p}2}\|_{L^2})^{\frac12-\frac3{2p}}
(N^{-\alpha_p}\|r^{\frac{p-3}2}|\tilde{V}_1|^{\frac{p}2}\|_{L^6})^{\frac7{2p}-\frac12}\\
&\times(N^{1-\alpha_p}\|r^{\frac{p-5}2}|\zeta_1|^{\frac{p}2}\|_{L^2})^{\frac12-\frac3{2p}}
(N^{-\alpha_p}\|r^{\frac{p-3}2}|\zeta_1|^{\frac{p}2}\|_{L^6})^{\frac7{2p}-\frac12}\\
&\times((2N)^{\beta_p}\|r^{\frac{p-3}2}|\zeta_2|^{\frac{p}2}\|_{L^2}+(2N)^{1+\beta_p}\|r^{\frac{p-5}2}|\zeta_2|^{\frac{p}2}\|_{L^2})\\
&\times((2N)^{\beta_p}\|r^{\frac{p-3}2}|\zeta_2|^{\frac{p}2}\|_{L^2})^{\frac2p}
((2N)^{\beta_p}\|r^{\frac{p-3}2}|\zeta_2|^{\frac{p}2}\|_{L^6})^{1-\frac4p}dt'\\
&\lesssim N^{-1+\frac3p+\frac2p\beta_p+\frac4p\alpha_p}\mathcal{H}_p^{1+\frac1p}(t).
\end{align*}
As for the estimate of the other part in $III_2^1$, namely
\begin{align*}
(2N)^{2\beta_p}\sum_{|k_1-k_2|=2}\int_0^t\int_{\mathbb{R}^3}|\tilde{V}_{k_1}||\zeta_{k_2}|(|\tilde{\nabla}\zeta_2|&+\frac1r|\zeta_2|)|\zeta_2|^{p-2}r^{p-3}dxdt'\\
&=(2N)^{2\beta_p}\sum_{j=1}^\infty\int_0^t\int_{\mathbb{R}^3}|\tilde{V}_j||\zeta_{j+2}|(|\tilde{\nabla}\zeta_2|+\frac1r|\zeta_2|)|\zeta_2|^{p-2}r^{p-3}dxdt'\\
&+(2N)^{2\beta_p}\sum_{j=1}^\infty\int_0^t\int_{\mathbb{R}^3}|\tilde{V}_{j+2}||\zeta_j|(|\tilde{\nabla}\zeta_2|+\frac1r|\zeta_2|)|\zeta_2|^{p-2}r^{p-3}dxdt',
\end{align*}
we observe that it has the same structure as that of $N^{-2\alpha_p}\int_0^tIII_1^1(t')dt'$. Then it follows
from the same derivation of $III_1^1$ that it can also be bounded by $N^{-1+\frac3p-\frac2p\beta_p}\mathcal{H}_p^{1+\frac1p}$. By summarizing the above estimates and using the fact that
\begin{eqnarray*}
-1+\frac3p+\frac2p\beta_p+\frac4p\alpha_p<0,
\end{eqnarray*}
which is guaranteed by the assumption $\alpha_p$, we achieve
\begin{align*}
(2N)^{2\beta_p}\int_0^tIII_2^1(t')dt'&\lesssim (N^{-1+\frac3p+\frac2p\beta_p+\frac4p\alpha_p}+N^{-1+\frac3p-\frac2p\beta_p})\mathcal{H}_p^{1+\frac1p}(t)\\
&\lesssim \mathcal{H}_p^{1+\frac1p}(t).
\end{align*}

\textbf{The estimate of $III_k^1$ when $k\geq3$.}

The case for $k\geq3$ is more complicated. And we shall decompose $III_k^1$ into three parts.

\textbf{Case 1:} The terms with $k_1$ or $k_2$ equals 1, which we denote as
\begin{align*}
III_k^{1,1}(t)&=\int_{\mathbb{R}^3}|\tilde{V}_1||\zeta_{k\pm1}|(|\tilde{\nabla}\zeta_k|+\frac1r|\zeta_k|)|\zeta_k|^{p-2}r^{p-3}dx
+\int_{\mathbb{R}^3}|\tilde{V}_{k\pm1}||\zeta_1|(|\tilde{\nabla}\zeta_k|+\frac1r|\zeta_k|)|\zeta_k|^{p-2}r^{p-3}dx\\
&:=III_k^{1,1,1}(t)+III_k^{1,1,2}(t).
\end{align*}
By using H\"{o}lder's inequality, we find
\begin{align*}
(kN)^{2\beta_p}III_k^{1,1,1}&\lesssim(N^{1-\alpha_p}\|r^{\frac{p-5}2}|\tilde{V}_1|^{\frac{p}2}\|_{L^2})^{\frac12-\frac3{2p}}
(N^{-\alpha_p}\|r^{\frac{p-3}2}|\tilde{V}_1|^{\frac{p}2}\|_{L^6})^{\frac7{2p}-\frac12}\\
&\times(((k\pm1)N)^{1+\beta_p}\|r^{\frac{p-5}2}|\zeta_{k\pm1}|^{\frac{p}2}\|_{L^2})^{\frac12-\frac3{2p}}
(((k\pm1)N)^{\beta_p}\|r^{\frac{p-3}2}|\zeta_{k\pm1}|^{\frac{p}2}\|_{L^6})^{\frac7{2p}-\frac12}\\
&\times((kN)^{\beta_p}\|r^{\frac{p-3}2}|\tilde{\nabla}\zeta_k||\zeta_k|^{\frac{p}2-1}\|_{L^2}+(kN)^{1+\beta_p}\|r^{\frac{p-5}2}|\zeta_k|^{\frac{p}2}\|_{L^2})\\
&\times((kN)^{\beta_p}\|r^{\frac{p-3}2}|\zeta_k|^{\frac{p}2}\|_{L^2})^{\frac2p}
((kN)^{\beta_p}\|r^{\frac{p-3}2}|\zeta_k|^{\frac{p}2}\|_{L^6})^{1-\frac4p}N^{-\frac12+\frac3{2p}+\frac2p\alpha_p}(kN)^{-\frac12+\frac3{2p}},
\end{align*}
where we used the fact that $\frac{k\pm1}k\in[\frac23,\frac43]$ for $k\geq3$. By integrating the above inequality over $[0,t]$ and then summing up the resulting inequalities for $k\geq3$, we obtain
\begin{align*}
\sum_{k\geq3}(kN)^{2\beta_p}\int_0^tIII_k^{1,1,1}(t')dt'&\lesssim (N^{2-2\alpha_p}\|r^{\frac{p-5}2}|\tilde{V}_1|^{\frac{p}2}\|_{L_t^2(L^2)}^2)^{\frac14-\frac3{4p}}
(N^{-2\alpha_p}\|r^{\frac{p-3}2}|\tilde{V}_1|^{\frac{p}2}\|_{L_t^2(L^6)}^2)^{\frac7{4p}-\frac14}\\
&\times(\sum_{k\geq3}((k\pm1)N)^{2+2\beta_p}\|r^{\frac{p-5}2}|\zeta_{k\pm1}|^{\frac{p}2}\|_{L_t^2(L^2)})^{\frac14-\frac3{4p}}(\sum_{k\geq3}((k\pm1)N)^{2\beta_p}\|r^{\frac{p-3}2}|\zeta_{k\pm1}|^{\frac{p}2}\|_{L_t^2(L^6)})^{\frac7{4p}-\frac14}\\
&\times(\sum_{k\geq3}(kN)^{2\beta_p}\|r^{\frac{p-3}2}|\tilde{\nabla}\zeta_k||\zeta_k|^{\frac{p}2-1}\|_{L_t^2(L^2)}^2+(kN)^{2+2\beta_p}\|r^{\frac{p-5}2}|\zeta_k|^{\frac{p}2}\|_{L_t^2(L^2)}^2)^{\frac12}\\
&\times(\sum_{k\geq3}(kN)^{2\beta_p}\|r^{\frac{p-3}2}|\zeta_k|^{\frac{p}2}\|_{L_t^\infty(L^2)})^{\frac1p}
(\sum_{k\geq3}(kN)^{2\beta_p}\|r^{\frac{p-3}2}|\zeta_k|^{\frac{p}2}\|_{L_t^2(L^6)})^{\frac12-\frac2p}\\
&\times\sup_{k\geq3}k^{-\frac12+\frac{3}{2p}}N^{-1+\frac3p+\frac2p\alpha_p}.
\end{align*}
Due to $-\frac12+\frac3{2p}<0$, we have $\sup_{k\geq3}k^{-\frac12+\frac{3}{2p}}=3^{-\frac12+\frac3{2p}}$. As a result, it comes out
\begin{eqnarray*}
\sum_{k\geq3}(kN)^{2\beta_p}\int_0^tIII_k^{1,1,1}(t')dt'\lesssim N^{-1+\frac3p+\frac2p\alpha_p}\mathcal{H}_p^{1+\frac1p}.
\end{eqnarray*}

As for $III_k^{1,1,2}(t)$, we also have
\begin{align*}
(kN)^{2\beta_p}III_k^{1,1,2}&\lesssim(((k\pm1)N)^{1+\beta_p}\|r^{\frac{p-5}2}|\tilde{V}_{k\pm 1}|^{\frac{p}2}\|_{L^2})^{\frac12-\frac3{2p}}
(((k\pm1)N)^{\beta_p}\|r^{\frac{p-3}2}|\tilde{V}_{k\pm1}|^{\frac{p}2}\|_{L^6})^{\frac7{2p}-\frac12}(N^{1-\alpha_p}\|r^{\frac{p-5}2}|\zeta_1|^{\frac{p}2}\|_{L^2})^{\frac12-\frac3{2p}}\\
&\times(N^{-\alpha_p}\|r^{\frac{p-3}2}|\zeta_1|^{\frac{p}2}\|_{L^6})^{\frac7{2p}-\frac12}((kN)^{\beta_p}\|r^{\frac{p-3}2}|\tilde{\nabla}\zeta_k||\zeta_k|^{\frac{p}2-1}\|_{L^2}+(kN)^{1+\beta_p}\|r^{\frac{p-5}2}|\zeta_k|^{\frac{p}2}\|_{L^2})\\
&\times((kN)^{\beta_p}\|r^{\frac{p-3}2}|\zeta_k|^{\frac{p}2}\|_{L^2})^{\frac2p}
((kN)^{\beta_p}\|r^{\frac{p-3}2}|\zeta_k|^{\frac{p}2}\|_{L^6})^{1-\frac4p}N^{-\frac12+\frac3{2p}+\frac2p\alpha_p}(kN)^{-\frac12+\frac3{2p}},
\end{align*}
where we used the fact that $\frac{k\pm1}k\in[\frac23,\frac43]$ for $k\geq3$. Hence, the similar work insures that
\begin{eqnarray*}
\sum_{k\geq3}(kN)^{2\beta_p}\int_0^tIII_k^{1,1,2}(t')dt'\lesssim N^{-1+\frac3p+\frac2p\alpha_p}\mathcal{H}_p^{1+\frac1p}.
\end{eqnarray*}

\textbf{Case 2:} The terms with $k_1$ or $k_2$ lying in $[2,k-2], k\geq4$, and we denote this part as
\begin{align*}
III_k^{1,2}(t)&=\sum_{j=2}^{k-2}\int_{\mathbb{R}^3}|\tilde{V}_j||\zeta_{k\pm j}|(|\tilde{\nabla}\zeta_k|+\frac1r|\zeta_k|)|\zeta_k|^{p-2}r^{p-3}dx+\sum_{j=2}^{k-2}\int_{\mathbb{R}^3}|\tilde{V}_{k\pm j}||\zeta_j|(|\tilde{\nabla}\zeta_k|+\frac1r|\zeta_k|)|\zeta_k|^{p-2}r^{p-3}dx\\
&=III_k^{1,2,1}(t)+III_k^{1,2,2}(t).
\end{align*}
By applying H\"{o}lder's inequality, we obtain
\begin{align*}
(kN)^{2\beta_p}III_k^{1,2,1}&\lesssim\sum_{j=2}^{k-2}((jN)^{1+\beta_p}\|r^{\frac{p-5}2}|\tilde{V}_j|^{\frac{p}2}\|_{L^2})^{\frac12-\frac3{2p}}
((jN)^{\beta_p}\|r^{\frac{p-3}2}|\tilde{V}_j|^{\frac{p}2}\|_{L^6})^{\frac7{2p}-\frac12}(((k\pm j)N)^{1+\beta_p}\|r^{\frac{p-5}2}|\zeta_{k\pm j}|^{\frac{p}2}\|_{L^2})^{\frac12-\frac3{2p}}\\
&\times
(((k\pm j)N)^{\beta_p}\|r^{\frac{p-3}2}|\zeta_{k\pm j}|^{\frac{p}2}\|_{L^6})^{\frac7{2p}-\frac12}((kN)^{\beta_p}\|r^{\frac{p-3}2}|\tilde{\nabla}\zeta_k||\zeta_k|^{\frac{p}2-1}\|_{L^2}+(kN)^{1+\beta_p}\|r^{\frac{p-5}2}|\zeta_k|^{\frac{p}2}\|_{L^2})\\
&\times((kN)^{\beta_p}\|r^{\frac{p-3}2}|\zeta_k|^{\frac{p}2}\|_{L^2})^{\frac2p}
((kN)^{\beta_p}\|r^{\frac{p-3}2}|\zeta_k|^{\frac{p}2}\|_{L^6})^{1-\frac4p}(jN)^{-\frac12+\frac3{2p}-\frac2p\beta_p}((k\pm j)N)^{-\frac12+\frac3{2p}-\frac2p\beta_p}(kN)^{\frac2p\beta_p}.
\end{align*}
And the terms containing $j$ in above inequality can be bounded by
\begin{align*}
&\left(\sum_{j=2}^{k-2}(jN)^{2+2\beta_p}\|r^{\frac{p-5}2}|\tilde{V}_j|^{\frac{p}2}\|_{L^2}^2\right)^{\frac14-\frac3{4p}}
\left(\sum_{j=2}^{k-2}(jN)^{2\beta_p}\|r^{\frac{p-3}2}|\tilde{V}_j|^{\frac{p}2}\|_{L^6}^2\right)^{\frac7{4p}-\frac14}\\
&\times\left(\sum_{j=2}^{k-2}((k\pm j)N)^{2+2\beta_p}\|r^{\frac{p-5}2}|\zeta_{k\pm j}|^{\frac{p}2}\|_{L^2}^2\right)^{\frac14-\frac3{4p}}
\left(\sum_{j=2}^{k-2}((k\pm j)N)^{2\beta_p}\|r^{\frac{p-3}2}|\zeta_{k\pm j}|^{\frac{p}2}\|_{L^6}^2\right)^{\frac7{4p}-\frac14}\\
&\times\left(\sum_{j=2}^{k-2}((j(k\pm j))N^2)^{(-\frac12+\frac3{2p}-\frac2p\beta_p)\times\frac{p}{p-2}}\right)^{1-\frac2p},
\end{align*}
Observing that $\beta_p<\frac{p-1}4$ ensures that
\begin{eqnarray*}
s=(-\frac12+\frac3{2p}-\frac2p\beta_p)\times\frac{p}{p-2}>-1,
\end{eqnarray*}
so that there holds
\begin{eqnarray*}
\lim_{k\rightarrow+\infty}\sum_{j=2}^{k-2}(\frac{j(k\pm j)}{k^2})^s\times\frac1k=\int_0^1\tau^s(1\pm \tau)^sd\tau<+\infty,
\end{eqnarray*}
which in particular implies for any $k\geq4$ that
\begin{eqnarray*}
\sum_{j=2}^{k-2}j^s(k\pm j)^s\lesssim k^{1+2s}.
\end{eqnarray*}
As a result, it comes out
\begin{eqnarray*}
\left(\sum_{j=2}^{k-2}((j(k\pm j))N^2)^{(-\frac12+\frac3{2p}-\frac2p\beta_p)\times\frac{p}{p-2}}\right)^{1-\frac2p}\lesssim k^{\frac1p-\frac4p\beta_p}N^{-1+\frac3p-\frac4p\beta_p}.
\end{eqnarray*}
Hence, we find
\begin{align*}
(kN)^{2\beta_p}\int_0^tIII_k^{1,2,1}(t')dt'&\lesssim k^{\frac1p-\frac4p\beta_p}N^{-1+\frac3p-\frac2p\beta_p}\mathcal{H}_p^{\frac2p}((kN)^{2\beta_p}\|r^{\frac{p-3}2}|\tilde{\nabla}\zeta_k||\zeta_k|^{\frac{p}2-1}\|_{L_t^2(L^2)}^2+(kN)^{2+2\beta_p}\|r^{\frac{p-5}2}|\zeta_k|^{\frac{p}2}\|_{L_t^2(L^2)}^2)^{\frac12}\\
&\times((kN)^{2\beta_p}\|r^{\frac{p-3}2}|\zeta_k|^{\frac{p}2}\|_{L_t^\infty(L^2)}^2)^{\frac1p}
((kN)^{2\beta_p}\|r^{\frac{p-3}2}|\zeta_k|^{\frac{p}2}\|_{L^6}^2)^{\frac12-\frac2p}.
\end{align*}
Our choice of $\beta_p>1$ ensures that
\begin{eqnarray*}
(\frac1p-\frac2p\beta_p)p<-1,
\end{eqnarray*}
so that we get, by summing up for $k\geq4$, that
\begin{align*}
\sum_{k\geq4}(kN)^{2\beta_p}\int_0^tIII_k^{1,2,1}(t')dt'&\lesssim(\sum_{k\geq4}k^{(\frac1p-\frac2p\beta_p)p})^{\frac1p}N^{-1+\frac3p-\frac2p\beta_p}\mathcal{H}_p^{\frac2p}\\
&\times(\sum_{k\geq4}(kN)^{2\beta_p}(\|r^{\frac{p-3}2}|\tilde{\nabla}\zeta_k||\zeta_k|^{\frac{p}2-1}\|_{L_t^2(L^2)}^2+(kN)^{2+2\beta_p}\|r^{\frac{p-5}2}|\zeta_k|^{\frac{p}2}\|_{L_t^2(L^2)}^2))^{\frac12}\\
&\times(\sum_{k\geq4}(kN)^{2\beta_p}\|r^{\frac{p-3}2}|\zeta_k|^{\frac{p}2}\|_{L_t^\infty(L^2)}^2)^{\frac1p}
(\sum_{k\geq4}(kN)^{2\beta_p}\|r^{\frac{p-3}2}|\zeta_k|^{\frac{p}2}\|_{L^6}^2)^{\frac12-\frac2p}\\
&\lesssim N^{-1+\frac3p-\frac2p\beta_p}\mathcal{H}_p^{1+\frac1p}.
\end{align*}

As for $III_k^{1,2,2}(t)$, we also have
\begin{align*}
(kN)^{2\beta_p}III_k^{1,2,2}&\lesssim\sum_{j=2}^{k-2}(((k\pm j)N)^{1+\beta_p}\|r^{\frac{p-5}2}|\tilde{V}_{k\pm j}|^{\frac{p}2}\|_{L^2})^{\frac12-\frac3{2p}}
(((k\pm j)N)^{\beta_p}\|r^{\frac{p-3}2}|\tilde{V}_{k\pm j}|^{\frac{p}2}\|_{L^6})^{\frac7{2p}-\frac12}\\
&\times((jN)^{1+\beta_p}\|r^{\frac{p-5}2}|\zeta_j |^{\frac{p}2}\|_{L^2})^{\frac12-\frac3{2p}}
((jN)^{\beta_p}\|r^{\frac{p-3}2}|\zeta_j|^{\frac{p}2}\|_{L^6})^{\frac7{2p}-\frac12}\\
&\times((kN)^{\beta_p}\|r^{\frac{p-3}2}|\tilde{\nabla}\zeta_k||\zeta_k|^{\frac{p}2-1}\|_{L^2}+(kN)^{1+\beta_p}\|r^{\frac{p-5}2}|\zeta_k|^{\frac{p}2}\|_{L^2})\\
&\times((kN)^{\beta_p}\|r^{\frac{p-3}2}|\zeta_k|^{\frac{p}2}\|_{L^2})^{\frac2p}
((kN)^{\beta_p}\|r^{\frac{p-3}2}|\zeta_k|^{\frac{p}2}\|_{L^6})^{1-\frac4p}\\
&\times (jN)^{-\frac12+\frac3{2p}-\frac2p\beta_p}((k\pm j)N)^{-\frac12+\frac3{2p}-\frac2p\beta_p}(kN)^{\frac2p\beta_p}.
\end{align*}
Thus, it refers to that
\begin{eqnarray*}
\sum_{k\geq4}(kN)^{2\beta_p}\int_0^tIII_k^{1,2,2}(t')dt'\lesssim N^{-1+\frac3p-\frac2p\beta_p}\mathcal{H}_p^{1+\frac1p}.
\end{eqnarray*}

\textbf{Cases 3:} The terms with both $k_1$ and $k_2$ being bigger than $k-1$, and we denote this part as
\begin{align*}
III_k^{1,3}(t)&=\sum_{j=k-1}^\infty\int_{\mathbb{R}^3}|\tilde{V}_j||\zeta_{k+ j}|(|\tilde{\nabla}\zeta_k|+\frac1r|\zeta_k|)|\zeta_k|^{p-2}r^{p-3}dx+\sum_{j=k-1}^\infty\int_{\mathbb{R}^3}|\tilde{V}_{k+ j}||\zeta_j|(|\tilde{\nabla}\zeta_k|+\frac1r|\zeta_k|)|\zeta_k|^{p-2}r^{p-3}dx\\
&=III_k^{1,3,1}(t)+III_k^{1,3,2}(t).
\end{align*}
Then we get
\begin{align*}
(kN)^{2\beta_p}III_k^{1,3,1}&\lesssim\sum_{j=k-1}^\infty((jN)^{1+\beta_p}\|r^{\frac{p-5}2}|\tilde{V}_j|^{\frac{p}2}\|_{L^2})^{\frac12-\frac3{2p}}
((jN)^{\beta_p}\|r^{\frac{p-3}2}|\tilde{V}_j|^{\frac{p}2}\|_{L^6})^{\frac7{2p}-\frac12}(((k+j)N)^{1+\beta_p}\|r^{\frac{p-5}2}|\zeta_{k+ j}|^{\frac{p}2}\|_{L^2})^{\frac12-\frac3{2p}}\\
&\times
(((k+j)N)^{\beta_p}\|r^{\frac{p-3}2}|\zeta_{k+ j}|^{\frac{p}2}\|_{L^6})^{\frac7{2p}-\frac12}((kN)^{\beta_p}\|r^{\frac{p-3}2}|\tilde{\nabla}\zeta_k||\zeta_k|^{\frac{p}2-1}\|_{L^2}+(kN)^{1+\beta_p}\|r^{\frac{p-5}2}|\zeta_k|^{\frac{p}2}\|_{L^2})\\
&\times((kN)^{\beta_p}\|r^{\frac{p-3}2}|\zeta_k|^{\frac{p}2}\|_{L^2})^{\frac2p}
((kN)^{\beta_p}\|r^{\frac{p-3}2}|\zeta_k|^{\frac{p}2}\|_{L^6})^{1-\frac4p}\\
&\times (jN)^{-\frac12+\frac3{2p}-\frac2p\beta_p}((k+j)N)^{-\frac12+\frac3{2p}-\frac2p\beta_p}(kN)^{\frac2p\beta_p}.
\end{align*}
Due to $\frac{k+j}j\in(1,3)$ for any $j\geq k-1$, ones have
\begin{align*}
(\sum_{j=k-1}^\infty (j(k+j))^{(-\frac12+\frac3{2p}-\frac2p\beta_p)\times\frac{p}{p-2}})^{1-\frac2p}&\lesssim
(\sum_{j=k-1}^\infty j^{(-1+\frac3p-\frac4p\beta_p)\times\frac{p}{p-2}})^{1-\frac2p}\\
&\lesssim k^{\frac1p-\frac4p\beta_p},
\end{align*}
from which, we deduce
\begin{align*}
(kN)^{2\beta_p}\int_0^tIII_k^{1,3,1}(t')dt'&\lesssim k^{\frac1p-\frac4p\beta_p}N^{-1+\frac3p-\frac2p\beta_p}\mathcal{H}_p^{\frac2p}((kN)^{2\beta_p}\|r^{\frac{p-3}2}|\tilde{\nabla}\zeta_k||\zeta_k|^{\frac{p}2-1}\|_{L_t^2(L^2)}^2+(kN)^{2+2\beta_p}\|r^{\frac{p-5}2}|\zeta_k|^{\frac{p}2}\|_{L_t^2(L^2)}^2)^{\frac12}\\
&\times((kN)^{2\beta_p}\|r^{\frac{p-3}2}|\zeta_k|^{\frac{p}2}\|_{L_t^\infty(L^2)}^2)^{\frac1p}
((kN)^{2\beta_p}\|r^{\frac{p-3}2}|\zeta_k|^{\frac{p}2}\|_{L^6}^2)^{\frac12-\frac2p}.
\end{align*}
By summing up the above inequality for $k\geq3$, we find
\begin{eqnarray*}
\sum_{k\geq3}(kN)^{2\beta_p}\int_0^tIII_k^{1,3,1}(t')dt'\lesssim N^{-1+\frac3p-\frac2p\beta_p}\mathcal{H}_p^{1+\frac1p}.
\end{eqnarray*}
At the same time, we also obtain
\begin{eqnarray*}
\sum_{k\geq3}(kN)^{2\beta_p}\int_0^tIII_k^{1,3,2}(t')dt'\lesssim N^{-1+\frac3p-\frac2p\beta_p}\mathcal{H}_p^{1+\frac1p}.
\end{eqnarray*}

Therefor, summarizing these estimations, we achieve
\begin{eqnarray}\label{eq17}
\int_0^tN^{-2\alpha_p}III_1^1+\sum_{k\geq2}(kN)^{2\beta_p}III_k^1dt'\lesssim \mathcal{H}_p^{1+\frac1p}.
\end{eqnarray}

\textbf{The estimate of $III_k^2$.}

Finally, let us handle the remaining term $III_k^2$. According to the relations between $k_1,k_2$ and $k$, we shall decompose into
\begin{eqnarray*}
III_k^2=III_k^{2,1}+III_k^{2,2},
\end{eqnarray*}
where
\begin{align*}
III_k^{2,1}=\sum_{(k_1,k_2)\in\Omega_k,~k_2\leq2k}\int_{\mathbb{R}^3}(k_1+k_2)N|U_{k_1}^\theta ||\zeta_{k_2}||\tilde{V}_k|^{p-1}r^{p-4}dx,~~III_k^{2,2}=\sum_{j=2k+1}\int_{\mathbb{R}^3}(k_1+k_2)N|U_{k_1}^\theta ||\zeta_{k_2}||\tilde{V}_k|^{p-1}r^{p-4}dx.
\end{align*}

Firstly, we focus on $III_k^{2,1}$. There holds for any $(k_1,k_2)\in\Omega_k$ with $k_2\leq2k$ that
\begin{eqnarray*}
k_1+k_2\leq k+2k_2\leq5k.
\end{eqnarray*}
Then it is easy to verify that $III_k^{2,1}$ shares a similar estimate as $III_k^1$. In fact, ones only needs to modify the estimate of $\tilde{V}_{k_1}$ in $III_k^1$ to $U_{k_1}^\theta$ in $III_k^2$, which gives an additional $N^{-\frac12}$ due to the different weight of $\tilde{V}_{k_1}$ and $U_{k_1}^\theta$ in the definition of $E_p$. Hence we have
\begin{eqnarray*}
\int_0^tN^{-2\alpha_p}III_1^{2,1}+\sum_{k\geq2}(kN)^{2\beta_p}III_k^{2,1}dt'\lesssim N^{-\frac12}\mathcal{H}_p^{1+\frac1p}.
\end{eqnarray*}

For $III_k^{2,2}$, we get, by using H\"{o}lder inequality, that
\begin{align*}
N^{-2\alpha_p}III_1^{2,2}&\lesssim\sum_{j=3}^\infty(((j\pm1)N)^{1+\frac{p}4+\beta_p}\|r^{\frac{p-5}2}|U_{j\pm1}^\theta|^{\frac{p}2}\|_{L^2})^{\frac2p}
((jN)^{1+\beta_p}\|r^{\frac{p-5}2}|\zeta_j|^{\frac{p}2}\|_{L^2})^{\frac2p}(N^{1-\alpha_p}\|r^{\frac{p-5}2}|\zeta_1|^{\frac{p}2}\|_{L^2})^{2-\frac7p}\\
&\times
(N^{-\alpha_p}\|r^{\frac{p-5}2}|\zeta_1|^{\frac{p}2}\|_{L^2})^{\frac2p}(N^{-\alpha_p}\|r^{\frac{p-3}2}|\zeta_1|^{\frac{p}2}\|_{L^6})^{\frac3p}j^{\frac12-\frac4p(1+\beta_p)}N^{-\frac32+\frac3p-\frac2p\alpha_p-\frac4p\beta_p},
\end{align*}
and for $k\geq2$ that
\begin{align*}
(kN)^{2\beta_p}III_k^{2,2}&\lesssim\sum_{j=3}^\infty(((j\pm k)N)^{1+\frac{p}4+\beta_p}\|r^{\frac{p-5}2}|U_{j\pm k}^\theta|^{\frac{p}2}\|_{L^2})^{\frac2p}
((jN)^{1+\beta_p}\|r^{\frac{p-5}2}|\zeta_j|^{\frac{p}2}\|_{L^2})^{\frac2p}\\
&\times((kN)^{1+\beta_p}\|r^{\frac{p-5}2}|\zeta_k|^{\frac{p}2}\|_{L^2})^{2-\frac7p}
((kN)^{\beta_p}\|r^{\frac{p-5}2}|\zeta_k|^{\frac{p}2}\|_{L^2})^{\frac2p}\\
&\times((kN)^{\beta_p}\|r^{\frac{p-3}2}|\zeta_k|^{\frac{p}2}\|_{L^6})^{\frac3p}j^{\frac12-\frac4p(1+\beta_p)}k^{-2+\frac7p+\frac2p\beta_p}N^{-\frac32+\frac3p-\frac2p\beta_p},
\end{align*}
where we used the fact that $2j\pm k\sim j\pm k\sim j$ for any $j\geq2k+1$. By integrating the above inequalities over $[0,t]$ and summing up the resulting inequalities for $k\geq1$, we obtain
\begin{align*}
\int_0^tN^{-2\alpha_p}III_1^{2,2}+(kN)^{2\beta_p}III_k^{2,2}dt'&\lesssim
\mathcal{H}_p^{1+\frac1p}N^{-\frac32+\frac3p-\frac2p\alpha_p-\frac4p\beta_p}(\sum_{j=3}^\infty j^{(\frac12-\frac4p(1+\beta_p))\times\frac{p}{p-2}})^{1-\frac2p}+\mathcal{H}_p^{\frac2p}N^{-\frac32+\frac3p-\frac2p\beta_p}\\
&\times\sum_{k\geq2}\{((kN)^{2+2\beta_p}\|r^{\frac{p-5}2}|\zeta_k|^{\frac{p}2}\|_{L^2}^2)^{1-\frac7{2p}}
((kN)^{2\beta_p}\|r^{\frac{p-5}2}|\zeta_k|^{\frac{p}2}\|_{L^2}^2)^{\frac1p}\\
&\times((kN)^{2\beta_p}\|r^{\frac{p-3}2}|\zeta_k|^{\frac{p}2}\|_{L^6}^2)^{\frac3{2p}}k^{-2+\frac7p+\frac2p\beta_p} (\sum_{j=2k+1}^\infty j^{(\frac12-\frac4p(1+\beta_p))\times\frac{p}{p-2}})^{1-2p}\}.
\end{align*}
It follows from $p<6$ and $\beta_p>1$ that
\begin{eqnarray*}
(\frac12-\frac4p(1+\beta_p))\times\frac{p}{p-2}<-1,
\end{eqnarray*}
so that one has
\begin{eqnarray*}
(\sum_{j=3}^\infty j^{(\frac12-\frac4p(1+\beta_p))\times\frac{p}{p-2}})^{1-\frac2p}\lesssim1,~~~
(\sum_{j=2k+1}^\infty j^{(\frac12-\frac4p(1+\beta_p))\times\frac{p}{p-2}})^{1-2p}\lesssim k^{\frac32-\frac4p(\frac32+\beta_p)}.
\end{eqnarray*}
Therefore, we obtain
\begin{align*}
\int_0^tN^{-2\alpha_p}III_1^{2,2}+\sum_{k\geq2}(kN)^{2\beta_p}III_k^{2,2}dt'&\lesssim N^{-\frac32+\frac3p-\frac2p\alpha_p-\frac4p\beta_p}\mathcal{H}_p^{1+\frac1p}
+N^{-\frac32+\frac3p-\frac2p\beta_p}\mathcal{H}_p^{1+\frac1p}(\sum_{k\geq2}k^{(-\frac12+\frac1p-\frac2p\beta_p)p})^{\frac1p}\\
&\lesssim N^{-\frac32+\frac3p-\frac2p\beta_p}\mathcal{H}_p^{1+\frac1p},
\end{align*}
where in the last step, we used the fact $(-\frac12+\frac1p-\frac2p\beta_p)p<-1$ for $\beta_p>1$ and $p>5$.

By summarizing these estimations, we acquire that
\begin{eqnarray}\label{eq18}
\int_0^tN^{-2\alpha_p}III_1^2+\sum_{k\geq2}(kN)^{2\beta_p}III_k^2dt'\lesssim \mathcal{H}_p^{1+\frac1p}.
\end{eqnarray}
Thanks to \eqref{eq16}-\eqref{eq18}, we complete the proof.$\square$

\section{The proof of Proposition 3.2}
The goal of this section is to present the proof of Proposition 3.2, which will be based on Proposition 3.1, Plancherel's identity and Hausdorff-Young's inequality.

\textbf{Proof of Proposition 3.2:} This process is made up of three steps. Firstly, we need the local well-posedness, which insures the existence of strong solutions. At last, prove the local solutions is global. In the Euclidean coordinates, the equations \eqref{eq4} becomes:
\begin{eqnarray*}
\begin{cases}
\partial_tu+u\cdot\nabla u-\Delta u+\nabla\pi=\frac{\eta e_3}{1+r^2},~~~&(t,x)\in(0,T]\times\mathbb{R}^3;\\
\partial_t\eta+u\cdot\nabla\eta-\Delta\eta-\frac{2x_hu^h\eta}{1+r^2}+\frac{4x_h\cdot\nabla_h\eta}{1+r^2}+\frac{4(1-r^2)\eta}{(1+r^2)^2}=0,~~~&(t,x)\in(0,T]\times\mathbb{R}^3;\\
div~u=0,~~~&(t,x)\in(0,T]\times\mathbb{R}^3;\\
u(0,x)=u_0(x),\eta(0,x)=\eta_0(x)~~~&x\in\mathbb{R}^3;
\end{cases}
\end{eqnarray*}
where $x_h=(x_1,x_2), \nabla_h=(\partial_1,\partial_2)$ and $u^h=(u^1,u^2), u=(u^h,u^3)$.

\textbf{Step 1: $L^2$ energy estimations}

By taking $L^2$ inner product with $u,\eta$ respectively, and using the divergence-free condition of $u$, we obtain
\begin{eqnarray*}
\frac12\frac{d}{dt}\|u(t)\|_{L^2}^2+\|\nabla u\|_{L^2}^2=\int_{\mathbb{R}^3}\frac{\eta u^3}{1+r^2}dx,
\end{eqnarray*}
and
\begin{eqnarray*}
\frac12\frac{d}{dt}\|\eta(t)\|_{L^2}^2+\|\nabla \eta\|_{L^2}^2=\int_{\mathbb{R}^3}(\frac{2x_hu^h\eta}{1+r^2}-\frac{4x_h\cdot\nabla_h\eta}{1+r^2}-\frac{4(1-r^2)\eta}{(1+r^2)^2})\eta dx.
\end{eqnarray*}
H\"{o}lder inequality gives that
\begin{align*}
|\int_{\mathbb{R}^3}(\frac{2(x_1u^1+x_2u^2)}{1+r^2}\eta^2dx|&\leq\int_{\mathbb{R}^3}|u^h||\eta|^2dx \leq\|u^h\|_{L^2}\|\eta\|_{L^3}\|\eta\|_{L^6}\leq C\| u^h\|_{L^2}\|\eta\|_{L^2}^{\frac12}\|\nabla\eta\|_{L^2}^{\frac32}\\
&\leq\frac14\|\nabla\eta\|_{L^2}^2+C\|u^h\|_{L^2}^4\|\eta\|_{L^2}^2,
\end{align*}
and
\begin{align*}
|-\int_{\mathbb{R}^3}(\frac{4x_h\cdot\nabla_h\eta}{1+r^2}+\frac{4(1-r^2)\eta}{(1+r^2)^2})\eta dx|&\leq2\int_{\mathbb{R}^3}|\nabla_h\eta||\eta|dx+4\|\eta\|_{L^2}^2\\
&\leq\frac14\|\nabla\eta\|_{L^2}^2+C\|\eta\|_{L^2}^2.
\end{align*}
Hence, we have
\begin{eqnarray}\label{eq19}
\frac{d}{dt}(\|\eta(t)\|_{L^2}^2+\|u(t)\|_{L^2}^2)+\|\nabla\eta\|_{L^2}^2+\|\nabla u\|_{L^2}^2\leq C(\|u\|_{L^2}^2+\|\eta\|_{L^2}^2)^3+C(\|u\|_{L^2}^2+\|\eta\|_{L^2}^2).
\end{eqnarray}

\textbf{Part two: $H^1$ energy estimations}

Nextly, we estimate the $H^1$ norm. By taking $L^2$ inner product with $-\Delta u,-\Delta\eta$ respectively, and using the divergence-free condition of $u$, we obtain
\begin{align*}
\frac12\frac{d}{dt}(\|\nabla u\|_{L^2}^2+\|\nabla \eta\|_{L^2}^2)&+\|\Delta u\|_{L^2}^2+\|\Delta \eta\|_{L^2}^2=-\int_{\mathbb{R}^3}(u\cdot\nabla u)\cdot\Delta u dx-\int_{\mathbb{R}^3}(u\cdot\nabla \eta)\cdot\Delta \eta dx\\
&+\int_{\mathbb{R}^3}\frac{\eta}{1+r^2}\cdot\Delta u^3 dx+\int_{\mathbb{R}^3}(\frac{2x_hu^h\eta}{1+r^2}-\frac{4x_h\cdot\nabla_h\eta}{1+r^2}-\frac{4(1-r^2)\eta}{1+r^2})\Delta\eta dx.
\end{align*}
H\"{o}lder inequality gives that
\begin{align*}
|\int_{\mathbb{R}^3}u\cdot\nabla u\cdot\Delta u dx+\int_{\mathbb{R}^3}u\cdot\nabla \eta\cdot\Delta \eta dx|&\lesssim\|u\|_{L^6}(\|\nabla u\|_{L^3}\|\Delta u\|_{L^2}+\|\nabla \eta\|_{L^3}\|\Delta \eta\|_{L^2})\\
&\lesssim\|\nabla u\|_{L^2}(\|\nabla u\|_{L^2}^{\frac12}\|\Delta u\|_{L^2}^{\frac32}+\|\nabla \eta\|_{L^2}^{\frac12}\|\Delta \eta\|_{L^2}^{\frac32})\\
&\leq\frac14\|\Delta u\|_{L^2}^2+\frac14\|\Delta \eta\|_{L^2}^2+C(\|\nabla u\|_{L^2}^2+\|\nabla \eta\|_{L^2}^2)^3,
\end{align*}
and
\begin{align*}
|\int_{\mathbb{R}^3}\frac{\eta}{1+r^2}\cdot\Delta u^3 dx&+\int_{\mathbb{R}^3}(\frac{2x_hu^h\eta}{1+r^2}-\frac{4x_h\cdot\nabla_h\eta}{1+r^2}-\frac{4(1-r^2)\eta}{1+r^2})\Delta\eta dx|\\
&\lesssim\|\eta\|_{L^2}\|\Delta u\|_{L^2}+\|\nabla\eta\|_{L^2}\|\Delta\eta\|_{L^2}+\|\eta\|_{L^2}\|\Delta\eta\|_{L^2}+\|u^h\|_{L^3}\|\eta\|_{L^6}\|\Delta\eta\|_{L^2}\\
&\lesssim\|\eta\|_{L^2}\|\Delta u\|_{L^2}+\|\nabla\eta\|_{L^2}\|\Delta\eta\|_{L^2}+\|\eta\|_{L^2}\|\Delta\eta\|_{L^2}+\|u^h\|_{L^2}^{\frac12}\|\nabla u^h\|_{L^2}^{\frac12}\|\nabla\eta\|_{L^2}\|\Delta\eta\|_{L^2}\\
&\leq\frac14\|\Delta u\|_{L^2}^2+\frac14\|\Delta \eta\|_{L^2}^2+C(1+\|u^h\|_{L^2}\|\nabla u^h\|_{L^2})\|\nabla\eta\|_{L^2}^2+C(\|u\|_{L^2}^2+\|\eta\|_{L^2}^2).
\end{align*}
Summing theses inequalities, we find that
\begin{align}\label{eq20}
\frac{d}{dt}(\|\nabla u\|_{L^2}^2+\|\nabla \eta\|_{L^2}^2)&+\|\Delta u\|_{L^2}^2+\|\Delta \eta\|_{L^2}^2\notag\\
&\leq C(\|\nabla u\|_{L^2}^2+\|\nabla \eta\|_{L^2}^2)^3+C(\|\nabla u\|_{L^2}^2+\|\nabla\eta\|_{L^2}^2)+C(\|u\|_{L^2}^2+\|\eta\|_{L^2}^2).
\end{align}
Combining \eqref{eq19} with \eqref{eq20}, ones conclude that
\begin{eqnarray*}
\frac{d}{dt}(\|u\|_{H^1}^2+\|\eta\|_{H^1}^2)+\|\nabla u\|_{H^1}^2+\|\nabla \eta\|_{H^1}^2\leq C(\|u\|_{H^1}^2+\|\eta\|_{H^1}^2)^3+C(\| u\|_{H^1}^2+\|\eta\|_{H^1}^2).
\end{eqnarray*}
Integrating them with respect to time over $[0,t]$, we have,
\begin{eqnarray*}
\|u(t)\|_{H^1}^2+\|\eta(t)\|_{H^1}^2\leq\|u_0\|_{H^1}^2+\|\eta_0\|_{H^1}^2+C\int_0^t(\| u\|_{H^1}^2+\|\eta\|_{H^1}^2)dt'+C\int_0^t(\| u\|_{H^1}^2+\|\eta\|_{H^1}^2)^3dt'.
\end{eqnarray*}
If time $T_l$ and initial data $u_0,\eta_0$ satisfy
\begin{eqnarray}\label{eq21}
CT_l[4(\|u_0\|_{H^1}^2+\|\eta_0\|_{H^1}^2)^2+1]\leq\frac12,
\end{eqnarray}
then for all $t\in(0,T_l]$,
\begin{eqnarray}\label{eq22}
\|u\|_{L_{T_l}^\infty(H^1)}^2+\|\nabla \eta\|_{L_{T_l}^\infty(H^1)}^2+\|\nabla u\|_{L_{T_l}^2(H^1)}^2+\|\nabla \eta\|_{L_{T_l}^2(H^1)}^2\leq2\|u_0\|_{H^1}^2+2\| \eta_0\|_{H^1}^2.
\end{eqnarray}
The estimation \eqref{eq22} suggest that the equations \eqref{eq4} with the initial data \eqref{eq12} has a unique solution $(u,\eta)$ on $[0,T_l]$, where $T_l$ is dependent of initial data. Furthermore, the solution can be written in terms of the Fourier series solutions \eqref{eq10}-\eqref{eq11}.

\textbf{Step three: Extension}

In Next, we show how to boost-strap the local analysis to any finite time $T$. Let $u_0(x)$ and $\eta_0$ be given by \eqref{eq12}, write
\begin{eqnarray*}
\alpha(T)=2(\|u_0\|_{L^2}^2+\|\eta_0\|_{L^2}^2)e^{C_0T},~~~\beta(T)=C_1\alpha(T)(1+T)
\end{eqnarray*}
where positive numbers $C_0$ and $C_1$ are only dependent of these embedding constants. For some $\epsilon>0$ to be determined later on, we denote the maximal time $T^*$ such that the following assumptions is right, namely that:
\begin{description}
 \item[i.] for all $0\leq t\leq T_l\leq T^*$,
  $$\|u\|_{L_{T_l}^\infty(L^2)}^2+\|\eta\|_{L_{T_l}^\infty(L^2)}^2+\|\nabla u\|_{L_{T_l}^2(L^2)}^2+\|\nabla\eta\|_{L_{T_l}^2(L^2)}^2\leq 2\alpha(T);$$
  \item[ii.] for all $0\leq t\leq T_l\leq T^*$,
  $$\|\nabla u\|_{L_{T_l}^\infty(L^2)}^2+\|\nabla \eta\|_{L_{T_l}^\infty(L^2)}^2+\|\Delta u\|_{L_{T_l}^2(L^2)}^2+\|\Delta \eta\|_{L_{T_l}^2(L^2)}^2\leq2\|\nabla u_0\|_{L^2}^2+2\|\nabla \eta_0\|_{L^2}^2+2\beta(T);$$
  \item[iii.] for all $0\leq t\leq T_l\leq T^*$, $\|U\|_{L_{T_l}^\infty(L^3)}+\|\xi\|_{L_{T_l}^\infty(L^3)}\leq2\epsilon$.
\end{description}
Since $U_0(x)=\xi_0(x)\equiv0$, $T^*$ is well-defined and must be positive. If $T^*\geq T$, then we achieve our purpose. Without loss of generality, we assume $T^*<T$ in what follows.

Denote $\bar{U}_k=(U_k^r,U_k^\theta,U_k^z)$ and $\bar{V}_k=(V_k^r,V_k^\theta,V_k^z)$. Applying the method in \cite{21} and the assumption \textbf{ii}, for $N$ large enough, we deduce that for any $t\leq T^*$
\begin{eqnarray}\label{eq23}
\sum_{k=1}^\infty\left((kN)^4\|(\frac{\bar{U}_k}{r^2},\frac{\bar{V}_k}{r^2})\|_{L_t^2(L^2)}^2+(kN)^2\|(\frac{\tilde{\nabla}\bar{U}_k}{r},\frac{\tilde{\nabla}\bar{V}_k}{r})\|_{L_t^2(L^2)}^2\right)\lesssim\|\nabla^2 u\|_{L_t^2(L^2)}^2\leq C_{in}(T)N^2.
\end{eqnarray}
where positive number $C_{in}(T)$ is only dependent of initial data and $T$. Similarly, the assumption \textbf{i} and Plancherel's identity that for all $t\leq T^*$
\begin{eqnarray}\label{eq24}
\sum_{k=1}^\infty(\|(\bar{U}_k,\bar{V}_k)\|_{L_t^\infty(L^2)}^2+\|(\tilde{\nabla}\bar{U}_k,\tilde{\nabla}\bar{V}_k)\|_{L_t^2(L^2)}^2+(kN)^2\|(\frac{\bar{U}_k}r,\frac{\bar{V}_k}r)\|_{L_t^2(L^2)}^2)\leq C_{in}(T).
\end{eqnarray}
By interpolating between \eqref{eq23} and \eqref{eq24}, we achieve for $t\leq T^*$
\begin{align}\label{eq25}
\sum_{k=1}^\infty k^3N^2\|(r^{-\frac32}\bar{U}_k,r^{-\frac32}\bar{V}_k)\|_{L_t^2(L^2)}^2
&+\sum_{k=1}^\infty k\|(r^{-\frac12}\tilde{\nabla}\bar{U}_k,r^{-\frac12}\tilde{\nabla}\bar{V}_k)\|_{L_t^2(L^2)}^2\notag\\
&\leq\left(\sum_{k=1}^\infty k^2N^2\|(r^{-1}\bar{U}_k,r^{-1}\bar{V}_k)\|_{L_t^2(L^2)}^2\right)^{\frac12}\left(\sum_{k=1}^\infty k^4N^2\|(r^{-2}\bar{U}_k,r^{-2}\bar{V}_k)\|_{L_t^2(L^2)}^2\right)^{\frac12}\notag\\
&+\left(\sum_{k=1}^\infty \|(\tilde{\nabla}\bar{U}_k,\tilde{\nabla}\bar{V}_k)\|_{L_t^2(L^2)}^2\right)^{\frac12}\left(\sum_{k=1}^\infty k^2\|(r^{-1}\tilde{\nabla}\bar{U}_k,r^{-1}\bar{V}_k)\|_{L_t^2(L^2)}^2\right)^{\frac12}\notag\\
&\leq C_{in}(T).
\end{align}

In view of the definition $T^*$, by taking $\epsilon$ to be so small that $\frac12-C(\|U\|_{L_{T^*}^\infty(L^3)}+\|\xi\|_{L_{T^*}^\infty(L^3)})\geq\frac13$, we deduce from Proposition 3.1 and \eqref{eq25} that for $t\leq T^*$, there holds
\begin{eqnarray}\label{eq26}
E_p(t)\leq3E_p(0)+CE_p^{1+\frac1p}+C_{in}(T)N^{-2}.
\end{eqnarray}
By the initial data, we find that for any $p\in(5,6)$
\begin{eqnarray*}
E_p(0)=N^{-2\alpha_p}\|r^{1-\frac3p}(a^r,b^r)\|_{L^p}^p+N^{-\frac{p}2-2\alpha_p}\|r^{1-\frac3p}(a^\theta,b^\theta)\|_{L^p}^p+N^{-2\alpha_p}\|r^{1-\frac3p}(a^z,b^z)\|_{L^p}^p
\leq C_{in}(T)N^{-2\alpha_p}.
\end{eqnarray*}
Then we deduce from \eqref{eq26} and a standard continuity argument that all $ t\leq T^*$,
\begin{eqnarray*}
E_p(t)\leq4E_p(0)+C_{in}(T)N^{-2}\leq C_{in}(T)N^{-2\alpha_p}.
\end{eqnarray*}
In particular, by taking $p=100/19,\alpha_p=1/30$ and $\beta_p=17/16$, and using the computations from \cite{21}, we obtain
\begin{align}\label{eq27}
&\|(\tilde{U}_1,\tilde{V}_1)\|_{L_t^5(L^5)}^5+\|(\xi_1,\zeta_1)\|_{L_t^5(L^5)}^5\leq C_{in}(T)N^{-2},~~~\|(U_1^\theta,V_1^\theta)\|_{L_t^5(L^5)}^5\leq C_{in}(T)N^{-2-75/31},\\
&\sum_{k=2}^\infty(kN)^4(\|(\tilde{U}_k,\tilde{V}_k)\|_{L_t^5(L^5)}^5+\|(\xi_k,\zeta_k)\|_{L_t^5(L^5)}^5)\leq C_{in}(T),~~~\sum_{k=2}^\infty(kN)^{4+75/31}\|(U_k^\theta,V_k^\theta)\|_{L_t^5(L^5)}^5\leq C_{in}(T),\notag
\end{align}
and
\begin{align}\label{eq28}
\|(\tilde{\nabla}\tilde{U}_1|\tilde{U}_1|^{\frac12},\tilde{\nabla}\tilde{V}_1|\tilde{V}_1|^{\frac12})\|_{L_t^2(L^2)}^2+\|(\tilde{\nabla}\xi_1|\xi_1|^{\frac12},\tilde{\nabla}\zeta_1|\zeta_1|^{\frac12})\|_{L_t^2(L^2)}^2&\leq C_{in}(T),\\
\sum_{k=2}^\infty k^{\frac43}(\|(\tilde{\nabla}\tilde{U}_k|\tilde{U}_k|^{\frac12},\tilde{\nabla}\tilde{V}_k|\tilde{V}_k|^{\frac12})\|_{L_t^2(L^2)}^2+\|(\tilde{\nabla}\xi_k|\xi_k|^{\frac12},\tilde{\nabla}\zeta_k|\zeta_k|^{\frac12})\|_{L_t^2(L^2)}^2)&\leq C_{in}(T)N^{-\frac23}.\notag
\end{align}

The rest of proof of Proposition 3.2 relies on the following three lemmas, which we admit for the time being.

\textbf{Lemma 5.1} There exists large enough integer $N_0$ so that for $N\geq N_0$ and for any $t\leq T^*$, there holds
\begin{eqnarray*}
\|U\|_{L_t^\infty(L^3)}^3+\||\nabla U||U|^{\frac12}\|_{L_t^2(L^2)}^2+\|\nabla |U|^{\frac32}\|_{L_t^2(L^2)}^2+\|\xi\|_{L_t^\infty(L^3)}^3+\||\nabla \xi||\xi|^{\frac12}\|_{L_t^2(L^2)}^2+\|\nabla |\xi|^{\frac32}\|_{L_t^2(L^2)}^2\leq4C_{in}(T)N^{-\frac35}.
\end{eqnarray*}
where the positive number $C_{in}(T)$ is only dependent of initial data and $T$.

\textbf{Lemma 5.2} There exists large enough integer $N_0$ so that for $N\geq N_0$ and for any $t\leq T^*$, there holds
\begin{eqnarray*}
\|u\|_{L_t^\infty(L^2)}^2+\|\eta\|_{L_t^\infty(L^2)}^2+\|\nabla u\|_{L_t^2(L^2)}^2+\|\nabla \eta\|_{L_t^2(L^2)}^2\leq\frac32\alpha(T).
\end{eqnarray*}

\textbf{Lemma 5.3} There exists large enough integer $N_0$ so that for $N\geq N_0$ and for any $t\leq T^*$, there holds
\begin{eqnarray*}
\|\nabla u\|_{L_t^\infty(L^2)}^2+\|\nabla \eta\|_{L_t^\infty(L^2)}^2+\|\Delta u\|_{L_t^2(L^2)}^2+\|\Delta \eta\|_{L_t^2(L^2)}^2\leq\frac32(\|\nabla u_0\|_{L^2}^2+\|\nabla \eta_0\|_{L^2}^2+\beta(T)).
\end{eqnarray*}

Now we are in a position to complete the proof. From the Lemma 5.1, for $N$ so large that $4C_{in}(T)N^{-\frac35}\leq\epsilon^3$, we deduce  that $\|U\|_{L_{T^*}^\infty(L^3)}+\|\xi\|_{L_{T^*}^\infty(L^3)}\leq\epsilon$, which together with Lemma 5.2 and Lemma 5.3 contradicts with the definition of $T^*$ unless $T^*\geq T$. This completes the proof. $\square$

Finally, we prove the three lemmas.

\textbf{Proof of Lemma 5.1:} Write
\begin{eqnarray*}
U(t,x)=U^r(t,r,z)e_r+U^\theta(t,r,z)e_\theta+U^z(t,r,z)e_z,~~~H(t,x)=H^r(t,r,z)e_r+H^\theta(t,r,z)e_\theta+H^z(t,r,z)e_z.
\end{eqnarray*}
Then in the Euclidean coordinates, the equations \eqref{eq6}-\eqref{eq7} becomes
\begin{eqnarray*}
\begin{cases}
\partial_tU+U\cdot\nabla U-\Delta U+\nabla\Pi=\frac{\xi e_3}{1+r^2}-H,~~~&(t,x)\in(0,T]\times\mathbb{R}^3;\\
\partial_t\xi+U\cdot\nabla\xi-\Delta\xi-\frac{2x_hU^h\xi}{1+r^2}+\frac{4x_h\cdot\nabla_h\xi}{1+r^2}+\frac{4(1-r^2)\xi}{(1+r^2)^2}=-\phi,~~~&(t,x)\in(0,T]\times\mathbb{R}^3;\\
div~U=0,~~~&(t,x)\in(0,T]\times\mathbb{R}^3;\\
U(0,x)=0,\xi(0,x)=0,~~~&x\in\mathbb{R}^3.
\end{cases}
\end{eqnarray*}
Let $\lambda>0$ and write for all $t\in[0,T^*]$
\begin{eqnarray*}
f_\lambda(t)=e^{-\lambda t}f(t),~~~f(t)=u(t),\eta(t),\Pi(t),H(t),\phi(t).
\end{eqnarray*}
By taking $L^2$ inner product with $|U_\lambda|U_\lambda$, and using the divergence-free condition of $U_\lambda$ and the fact that the initial data of $U$ vanishes, we obtain
\begin{align*}
\frac13\frac{d}{dt}(\|U_\lambda\|_{L^3}^3&+\|\xi_\lambda\|_{L^3}^3)+(\||\nabla U_\lambda||U_\lambda|^{\frac12}\|_{L^2}^2+\frac49\|\nabla |U_\lambda|^{\frac32}\|_{L^2}^2)+(\||\nabla \xi_\lambda||\xi_\lambda|^{\frac12}\|_{L^2}^2+\frac49\|\nabla |\xi_\lambda|^{\frac32}\|_{L^2}^2)\\
&\leq-\lambda(\|U_\lambda\|_{L^3}^3+\|\xi_\lambda\|_{L^3}^3)+\int_{\mathbb{R}^3}|\tilde{\nabla}\Pi_\lambda||U_\lambda|^2dx
+\|\xi_\lambda\|_{L^3}\|U_\lambda\|_{L^3}^2-\int_{\mathbb{R}^3}H_\lambda U_\lambda|U_\lambda|dx\\
&-\int_{\mathbb{R}^3}\phi_\lambda \xi_\lambda|\xi_\lambda|dx+\int_{\mathbb{R}^3}|U^h||\xi_\lambda|^3dx+2\int_{\mathbb{R}^3}|\nabla_h\xi_\lambda||\xi_\lambda|^2dx
+4\|\xi_\lambda\|_{L^3}^3.
\end{align*}
Nextly, we will estimate these terms one by one.

First of all, we determine the positive number $\lambda$. On the one hand, Holder inequality gives
\begin{align*}
\|\xi_\lambda\|_{L^3}\|U_\lambda\|_{L^3}^2+\int_{\mathbb{R}^3}|\nabla_h\xi_\lambda||\xi_\lambda|^2dx&\leq\||\nabla \xi_\lambda||\xi_\lambda|^{\frac12}\|_{L^2}\||\xi_\lambda|^{\frac32}\|_{L^2}+\|\xi_\lambda\|_{L^3}\|U_\lambda\|_{L^3}^2\\
&\leq\frac18\||\nabla \xi_\lambda||\xi_\lambda|^{\frac12}\|_{L^2}^2+C(\|\xi_\lambda\|_{L^3}^3+\|U_\lambda\|_{L^3}^3).
\end{align*}
On the another hands, we deal with the pressure term. Thanks to $div~U=0$, the pressure term $\Pi_\lambda$ meets
\begin{eqnarray*}
-\Delta\Pi_\lambda=e^{\lambda t}div(U_\lambda\cdot\nabla U_\lambda)+div~H_\lambda-\frac1{1+r^2}\partial_3\xi_\lambda.
\end{eqnarray*}
Hence, the integration including the pressure term $\Pi_\kappa$ is decomposed as
\begin{eqnarray*}
\int_{\mathbb{R}^3}\nabla\Pi_\kappa\cdot U_\kappa|U_\kappa|dx=\sum_{j=1}^2II_j(t),
\end{eqnarray*}
where
\begin{align*}
II_1(t)&=\int_{\mathbb{R}^3}\nabla(-\Delta)^{-1}(div(e^{\lambda t}U_\lambda\cdot\nabla U_\lambda)+\frac1{1+r^2}\partial_3\xi_\lambda)\cdot U_\lambda|U_\lambda|dx\\
II_2(t)&=\int_{\mathbb{R}^3}\nabla(-\Delta)^{-1}(div~H_\lambda)\cdot U_\lambda|U_\lambda|dx.
\end{align*}
The $L^p$ boundedness of the operators $\nabla(-\Delta)^{-1}\nabla$, we deduce that
\begin{align*}
|II_1|&\lesssim e^{\lambda t}\|U_\lambda\cdot\nabla U_\lambda\|_{L^{\frac32}}\|U_\lambda\|_{L^3}^{\frac12}\|U_\lambda\|_{L^9}^{\frac32}+\|\xi_\lambda\|_{L^3}\|U_\lambda\|_{L^3}^2\\
&\leq Ce^{\lambda t}\|U_\lambda\cdot\nabla U_\lambda\|_{L^{\frac32}}\|U_\lambda\|_{L^3}^{\frac12}\|U_\lambda\|_{L^9}^{\frac32}+C(\|U_\lambda\|_{L^3}^3+\|\xi_\lambda\|_{L^3}^3).
\end{align*}
where we use the $L^p$ boundedness of the operator $(-\Delta)^{-1}\nabla^2$, and
\begin{eqnarray*}
|II_2|\lesssim\|H\|_{L^{\frac53}}\|U_\lambda\|_{L^3}^{\frac45}\|U_\lambda\|_{L^9}^{\frac65}. \end{eqnarray*}
Thus, there is a positive number $\lambda_0$ such that when $\lambda\geq\lambda_0$, these terms $\|U_\lambda\|_{L^3}^3+\|\xi_\lambda\|_{L^3}^3$ can be absorbed.

Now given $\lambda\geq\lambda_0+8C$, due to $T^*<T$, it implies for all $T^*\geq t>0$, $e^{\lambda T^*}<e^{\lambda T}$ and we also have that
\begin{align*}
\frac13(\|U_\lambda\|_{L_t^\infty(L^3)}^3+\|\xi_\lambda\|_{L_t^\infty(L^3)}^3)&+(\||\nabla U_\lambda||U_\lambda|^{\frac12}\|_{L_t^2(L^2)}^2+\frac49\|\nabla |U_\lambda|^{\frac32}\|_{L_t^2(L^2)}^2)\\
&~~~~~~~+(\frac78\||\nabla\xi_\lambda||\xi_\lambda|^{\frac12}\|_{L_t^2(L^2)}^2+\frac49\|\nabla|\xi_\lambda|^{\frac32}\|_{L_t^2(L^2)}^2)+(\lambda-\lambda_0)(\|U_\lambda\|_{L_t^3(L^3)}^3+\|\xi_\lambda\|_{L_t^3(L^3)}^3)\\
&\leq C\int_0^t\|H\|_{L^{\frac53}}\|U_\lambda\|_{L^3}^{\frac45}\|U_\lambda\|_{L^9}^{\frac65}dt'
+Ce^{\lambda T}\int_0^t\|U_\lambda\cdot\nabla U_\lambda\|_{L^{\frac32}}\|U_\lambda\|_{L^3}^{\frac12}\|U_\lambda\|_{L^9}^{\frac32}dt'\\
&-\int_0^t\int_{\mathbb{R}^3}\phi_\lambda\xi_\lambda|\xi_\lambda|dxdt'+\int_0^t\int_{\mathbb{R}^3}|U^h||\xi_\lambda|^3dxdt'.
\end{align*}
It is easied to compute that
\begin{align*}
\int_0^t\int_{\mathbb{R}^3}|U^h||\xi_\lambda|^3dxdt'&\leq\int_0^t\|U^h\|_{L^3}\||\xi_\lambda|^{\frac32}\|_{L^6}\||\xi_\lambda|^{\frac32}\|_{L^2}dt'\\
&\leq \frac18\|U\|_{L_t^\infty(L^3)}\|\nabla|\xi_\lambda|^{\frac32}\|_{L^2(L^2)}^2+8C\|U\|_{L_t^\infty(L^3)}\|\xi_\lambda\|_{L_t^3(L^3)}^3\\
&\leq \frac18e^{\lambda T}\|U_\lambda\|_{L_t^\infty(L^3)}\|\nabla|\xi_\lambda|^{\frac32}\|_{L^2(L^2)}^2+8C\|U\|_{L_t^\infty(L^3)}\|\xi_\lambda\|_{L_t^3(L^3)}^3.
\end{align*}
and
\begin{align*}
\int_0^t\|U_\lambda\cdot\nabla U_\lambda\|_{L^{\frac32}}\|U_\lambda\|_{L^3}^{\frac12}\|U_\lambda\|_{L^9}^{\frac32}dt'
&\leq \|U_\lambda\cdot\nabla U_\lambda\|_{L_t^2(L^{\frac32})}\|U_\lambda\|_{L_t^\infty(L^3)}^{\frac12}\|U_\lambda\|_{L_t^3(L^9)}^{\frac32}\\
&\leq \|U_\lambda\|_{L_t^\infty(L^3)}\||U_\lambda|^{\frac12}|\nabla U_\lambda|\|_{L_t^2(L^2)}\|\nabla|U_\lambda|^{\frac32}\|_{L_t^2(L^2)}.
\end{align*}
The expression of vector $H$ and \eqref{eq27}-\eqref{eq28} shows that
\begin{align*}
\|H\|_{L_t^{\frac32}(L^{\frac53})}&\leq\sum_{k=1}^\infty(\||U_k|^{\frac12}\tilde{\nabla}U_k\|_{L_t^2(L^2)}\|U_k\|_{L_t^5(L^5)}^{\frac12}+\||V_k|^{\frac12}\tilde{\nabla}V_k\|_{L_t^2(L^2)}\|\tilde{V}_k\|_{L_t^5(L^5)}^{\frac12}\\
&+kN\|r^{-\frac32}U_k^\theta\|_{L_t^2(L^2)}^{\frac23}\|\tilde{V}_k\|_{L_t^5(L^5)}^{\frac13}\|U_k^\theta\|_{L_t^5(L^5)}
+kN\|r^{-\frac32}V_k^\theta\|_{L_t^2(L^2)}^{\frac23}\|\tilde{U}_k\|_{L_t^5(L^5)}^{\frac13}\|V_k^\theta\|_{L_t^5(L^5)}\\
&+\|r^{-\frac32}U_k^\theta\|_{L_t^2(L^2)}^{\frac23}\|U_k\|_{L_t^5(L^5)}^{\frac43}+\|r^{-\frac32}V_k^\theta\|_{L_t^2(L^2)}^{\frac23}\|V_k\|_{L_t^5(L^5)}^{\frac43})\\
&\leq C_{in}(T)N^{-\frac15}.
\end{align*}
Then
\begin{align*}
C\int_0^t\|H\|_{L^{\frac53}}\|U_\lambda\|_{L^3}^{\frac45}\|U_\lambda\|_{L^9}^{\frac65}dt'&\leq \|H\|_{L_t^{\frac32}(L^{\frac53})}\|U_\lambda\|_{L^3}^{\frac45}\|\nabla|U_\lambda|^{\frac32}\|_{L^2}^{\frac45}\\
&\leq C_{in}N^{-\frac15}\|U_\lambda\|_{L^3}^{\frac45}\|\nabla|U_\lambda|^{\frac32}\|_{L^2}^{\frac45}\\
&\leq\frac14(\frac13\|U_\lambda\|_{L_t^\infty(L^3)}^3+\frac49\|\nabla |U_\lambda|^{\frac32}\|_{L_t^2(L^2)}^2)+C_{in}(T)N^{-\frac35}.
\end{align*}
Finally, we estimate the rest term $\phi_\lambda\xi_\lambda|\xi_\lambda|$. Thanks to that
\begin{eqnarray*}
\int_{\mathbb{R}^3}(\tilde{U}_k\cdot\tilde{\nabla}\xi_k)\xi_\lambda|\xi_\lambda|dx =-\int_{\mathbb{R}^3}\frac{kN}rV_k^\theta\xi_k\xi_\lambda|\xi_\kappa|dx-\int_{\mathbb{R}^3}\xi_k\tilde{U}_k\cdot\tilde{\nabla}(\xi_\lambda|\xi_\lambda|)dx,
\end{eqnarray*}
thus, ones have
\begin{align*}
|\int_0^t\int_{\mathbb{R}^3}\phi\xi_\lambda|\xi_\lambda|dxdt'|&\leq 2\sum_{k=1}^\infty\int_0^t\int_{\mathbb{R}^3}(\frac{kN}r|U_k^\theta| |\zeta_k|+\frac{kN}r|V_k^\theta||\xi_k|)|\xi_\lambda|^2dxdt'\\
&+2\sum_{k=1}^\infty\int_0^t\int_{\mathbb{R}^3}|\xi_k||\tilde{U}_k||\tilde{\nabla}\xi_\lambda||\xi_\lambda|dxdt'\\
&:=III_1(t)+III_2(t).
\end{align*}
For the term $III_1(t)$, we can use the way to deal with $H$. Hence,
\begin{align*}
III_1(t)&=2\sum_{k=1}^\infty(kN\|r^{-\frac32}U_k^\theta\|_{L_t^2(L^2)}^{\frac23}\|\zeta_k\|_{L_t^5(L^5)}^{\frac13}\|U_k^\theta\|_{L_t^5(L^5)}+kN\|r^{-\frac32}V_k^\theta\|_{L_t^2(L^2)}^{\frac23}\|\xi_k\|_{L_t^5(L^5)}^{\frac13}\|V_k^\theta\|_{L_t^5(L^5)})\|\xi_\lambda\|_{L_t^\infty(L^3)}^{\frac45}\|\xi_\lambda\|_{L_t^3(L^9)}^{\frac65}\\
&\leq \frac14(\frac13\|\xi_\lambda\|_{L_t^\infty(L^3)}^3+\frac49\|\nabla |\xi_\lambda|^{\frac32}\|_{L_t^2(L^2)}^2)+C_{in}(T)N^{-\frac35}.
\end{align*}
Compared $III_1(t)$, the term $III_2(t)$ is more skillful, which is estimated by
\begin{align*}
III_2(t)&\leq C\sum_{k=1}^\infty\|\xi_k\|_{L_t^5(L^5)}\|\tilde{U}_k\|_{L_t^5(L^5)}\|\xi_\lambda\|_{L_t^\infty(L^3)}^{\frac15}\|\nabla|\xi_\lambda|^\frac32\|_{L_t^2(L^2)}^{\frac15}\||\xi_\lambda|^{\frac12}|\nabla\xi_\lambda|\|_{L_t^2(L^2)}\\
&\leq C_{in}(T)N^{-\frac45}\|\xi_\lambda\|_{L_t^\infty(L^3)}^{\frac15}\|\nabla|\xi_\lambda|^\frac32\|_{L_t^2(L^2)}^{\frac15}\||\xi_\lambda|^{\frac12}|\nabla\xi_\lambda|\|_{L_t^2(L^2)}\\
&\leq\frac14(\frac13\|\xi_\lambda\|_{L_t^\infty(L^3)}^3+\||\nabla\xi_\lambda| |\xi_\lambda|^{\frac12}\|_{L_t^2(L^2)}^2+\frac49\|\nabla |\xi_\lambda|^{\frac32}\|_{L_t^2(L^2)}^2)+C_{in}(T)N^{-\frac35}.
\end{align*}
Here, we use the fact
\begin{eqnarray*}
\|\xi_\lambda\|_{L_t^5(L^5)}^5\leq\int_0^t\|\xi_\lambda\|_{L^3}\||\xi_\lambda|^\frac32\|_{L^4}^{\frac83}dt'\leq\int_0^t\|\xi_\lambda\|_{L^3}\||\xi_\lambda|^\frac32\|_{L^2}^{\frac23}\||\xi_\lambda|^\frac32\|_{L^6}^2dt'
\leq C\|\xi_\lambda\|_{L_t^\infty(L^3)}^2\|\nabla|\xi_\lambda|^\frac32\|_{L_t^2(L^2)}^2.
\end{eqnarray*}

Therefore, summing up the above all estimations, we finally have
\begin{align*}
(\frac12-Ce^{\lambda T}\|U_\lambda\|_{L_t^\infty(L^3)})&(\frac13\|U_\lambda\|_{L_t^\infty(L^3)}^3+\||\nabla U_\lambda||U_\lambda|^{\frac12}\|_{L_t^2(L^2)}^2+\frac49\|\nabla |U_\lambda|^{\frac32}\|_{L_t^2(L^2)}^2)\\
&+\frac12(\frac13\|\xi_\lambda\|_{L_t^\infty(L^3)}^3+\||\nabla \xi_\lambda||\xi_\lambda|^{\frac12}\|_{L_t^2(L^2)}^2+\frac49\|\nabla |\xi_\lambda|^{\frac32}\|_{L_t^2(L^2)}^2)\\
&+C(1-\|U\|_{L_t^\infty(L^3)})(\|U_\lambda\|_{L_t^3(L^3)}^3+\|\xi_\lambda\|_{L_t^3(L^3)}^3)
\leq C_{in}(T)N^{-\frac35}.
\end{align*}
So ones acquire that
\begin{eqnarray*}
(\frac12-Ce^{\lambda T}\|U_\lambda\|_{L_t^\infty(L^3)})(\|U_\lambda\|_{L_t^\infty(L^3)}+\|\xi_\lambda\|_{L_t^\infty(L^3)})\leq C_{in}(T)N^{-\frac15},
\end{eqnarray*}
as long as $\|U\|_{L_t^\infty(L^3)}\leq1$ for all $t\in(0,T^*]$. Then by virtue of definition $T^*$ as long as $\epsilon$ is so small that $C\epsilon e^{\lambda T}\leq\frac14$, we achieve the purpose. $\square$

\textbf{Proof of Lemma 5.2:} We first get, by using Hausdorff-Young's inequality, that
\begin{eqnarray*}
\|u-U\|_{L^5((0,2\pi),d\theta)}\leq C\|(U_k,V_k)_{k\in\mathbb{N}^+}\|_{l^{\frac54}(\mathbb{N}^+)}.
\end{eqnarray*}
By taking $L_t^5(L^5(rdrdz))$ norm of the above inequality and then using Minkowski's inequality and \eqref{eq27}, we find
\begin{align}\label{eq29}
\|u-U\|_{L_t^5(L^5)}&\leq C\|(\|(U_k,V_k)\|_{L_t^5(L^5)})_{k\in\mathbb{N}^+}\|_{l^{\frac54}(\mathbb{N}^+)}\notag\\
&\leq C\|(\|k^{\frac45}(U_k,V_k)\|_{L_t^5(L^5)})_{k\in\mathbb{N}^+}\|_{l^5(\mathbb{N}^+)}\|(k^{-\frac45})_{k\in\mathbb{N}^+}\|_{l^{\frac53}(\mathbb{N}^+)}\notag\\
&\leq C_{in}(T)N^{-\frac25}.
\end{align}
While it follows Lemma 5.1, that
\begin{eqnarray*}
\|U\|_{L_t^5(L^5)}\leq\|U\|_{L_t^\infty(L^3)}^{\frac25}\|U\|_{L_t^3(L^9)}^{\frac35}\leq C\|U\|_{L_t^\infty(L^3)}^{\frac25}\|\nabla|U|^{\frac32}\|_{L_t^2(L^2)}^{\frac35}\leq C_{in}(T)N^{-\frac15}.
\end{eqnarray*}
which together with \eqref{eq29} ensures that for all $0<t\leq T^\ast$
\begin{eqnarray}\label{eq30}
\|u\|_{L_t^5(L^5)}\leq C_{in}(T)N^{-\frac15}.
\end{eqnarray}

By taking $L^2$ inner product with $u,\eta$ respectively, and using the divergence-free condition of $u$, we obtain
\begin{align*}
\frac12&\frac{d}{dt}(\|u\|_{L^2}^2+\|\eta\|_{L^2}^2)+\|\nabla u\|_{L^2}^2+\|\nabla\eta\|_{L^2}^2=
\int_{\mathbb{R}^3}\frac{u^3\eta}{1+r^2} dx+\int_{\mathbb{R}^3}(\frac{2x_hu^h\eta}{1+r^2}-\frac{4x_h\cdot\nabla_h\eta}{1+r^2}-\frac{4(1-r^2)\eta}{1+r^2})\eta dx.
\end{align*}
Against, H\"{o}lder inequality gives
\begin{align*}
|\int_{\mathbb{R}^3}\frac{u^3\eta}{1+r^2} dx&+\int_{\mathbb{R}^3}(\frac{2x_hu^h\eta}{1+r^2}-\frac{4x_h\cdot\nabla_h\eta}{1+r^2}-\frac{4(1-r^2)\eta}{1+r^2})\eta dx|\\
&\lesssim\|\eta\|_{L^2}\|u\|_{L^2}+\|\nabla\eta\|_{L^2}\|\eta\|_{L^2}+\|\eta\|_{L^2}^2+\|u\|_{L^3}\|\eta\|_{L^6}\|\eta\|_{L^2}\\
&\lesssim\|\eta\|_{L^2}\|u\|_{L^2}+\|\nabla\eta\|_{L^2}\|\eta\|_{L^2}+\|\eta\|_{L^2}^2+\|u\|_{L^3}\|\eta\|_{L^6}\|\eta\|_{L^2}+\|u\|_{L^2}^{\frac49}\| u\|_{L^5}^{\frac59}\|\nabla\eta\|_{L^2}\|\Delta\eta\|_{L^2}\\
&\leq\frac14\|\Delta \eta\|_{L^2}^2+C\|u\|_{L^2}^{\frac89}\|u\|_{L^5}^{\frac{10}9}\|\eta\|_{L^2}^2+C(\|u\|_{L^2}^2+\|\eta\|_{L^2}^2).
\end{align*}
Summing up these estimations, we find that
\begin{eqnarray*}
\frac{d}{dt}(\|u\|_{L^2}^2+\|\eta\|_{L^2}^2)+\|\nabla u\|_{L^2}^2+\|\nabla\eta\|_{L^2}^2\leq C\|u\|_{L^2}^{\frac89}\|u\|_{L^5}^{\frac{10}9}(\|u\|_{L^2}^2+\|\eta\|_{L^2}^2)+C(\|u\|_{L^2}^2+\|\eta\|_{L^2}^2).
\end{eqnarray*}
By applying Gronwall's inequality and using the estimation \eqref{eq30}, we achieve for all $t\leq T^*$
\begin{align*}
\|u\|_{L_t^\infty(L^2)}^2+\|\eta\|_{L_t^\infty(L^2)}^2+\|\nabla u\|_{L_t^2(L^2)}^2+\|\nabla \eta\|_{L_t^2(L^2)}^2&\leq(\|u_0\|_{L^2}^2+\| \eta_0\|_{L^2}^2)\exp\left(\int_0^tC_0+C\|u\|_{L^2}^{\frac89}\|u\|_{L^5}^{\frac{10}9}dt'\right)\\
&\leq(\|u_0\|_{L^2}^2+\|\eta_0\|_{L^2}^2)e^{C_0T}e^{C_{in}(T)N^{-2/9}}\\
&\leq\frac32\alpha(T),
\end{align*}
provided that $N$ is sufficiently large, where we use the fact:
\begin{eqnarray*}
\int_0^t\|u\|_{L^2}^{\frac89}\|u\|_{L^5}^{\frac{10}9}dt'\leq (T^*)^{\frac79}\|u\|_{L_t^\infty(L^2)}^{\frac89}\|u\|_{L_t^5(L^5)}^{\frac{10}9}\leq C_{in}(T)N^{-\frac29}.
\end{eqnarray*}
$\square$

\textbf{Proof of Lemma 5.3:} By taking $L^2$ inner product with $-\Delta u,-\Delta\eta$ respectively, and using the divergence-free condition of $u$, we obtain
\begin{align*}
\frac12\frac{d}{dt}(\|\nabla u\|_{L^2}^2+\|\nabla \eta\|_{L^2}^2)&+\|\Delta u\|_{L^2}^2+\|\Delta \eta\|_{L^2}^2=-\int_{\mathbb{R}^3}(u\cdot\nabla u)\cdot\Delta u dx-\int_{\mathbb{R}^3}(u\cdot\nabla \eta)\cdot\Delta \eta dx\\
&+\int_{\mathbb{R}^3}\frac{\eta}{1+r^2}\cdot\Delta u^3 dx+\int_{\mathbb{R}^3}(\frac{2x_hu^h\eta}{1+r^2}-\frac{4x_h\cdot\nabla_h\eta}{1+r^2}-\frac{4(1-r^2)\eta}{1+r^2})\Delta\eta dx.
\end{align*}
Against, H\"{o}lder inequality gives
\begin{align*}
|\int_{\mathbb{R}^3}(u\cdot\nabla u)\cdot\Delta u dx|+|\int_{\mathbb{R}^3}(u\cdot\nabla \eta)\cdot\Delta \eta dx|
&\leq C\|u\|_{L^5}\|\nabla u\|_{L^{\frac{10}3}}\|\Delta u\|_{L^2}+C\|u\|_{L^5}\|\nabla \eta\|_{L^{\frac{10}3}}\|\Delta\eta\|_{L^2}\\
&\leq C\|u\|_{L^5}\|\nabla u\|_{L^2}^{\frac25}\|\Delta u\|_{L^2}^{\frac85}+C\|u\|_{L^5}\|\nabla \eta\|_{L^2}^{\frac25}\|\Delta\eta\|_{L^2}^{\frac85}\\
&\leq\frac14\|\Delta u\|_{L^2}^2+\frac14\|\Delta \eta\|_{L^2}^2+C\|u\|_{L^5}^5(\|\nabla u\|_{L^2}^2+\|\nabla \eta\|_{L^2}^2),
\end{align*}
and
\begin{align*}
|\int_{\mathbb{R}^3}\frac{\eta}{1+r^2}\cdot\Delta u^3 dx&+\int_{\mathbb{R}^3}(\frac{2x_hu^h\eta}{1+r^2}-\frac{4x_h\cdot\nabla_h\eta}{1+r^2}-\frac{4(1-r^2)\eta}{1+r^2})\Delta\eta dx|\\
&\lesssim\|\eta\|_{L^2}\|\Delta u\|_{L^2}+\|\nabla\eta\|_{L^2}\|\Delta\eta\|_{L^2}+\|\eta\|_{L^2}\|\Delta\eta\|_{L^2}+\|u\|_{L^3}\|\eta\|_{L^6}\|\Delta\eta\|_{L^2}\\
&\lesssim\|\eta\|_{L^2}\|\Delta u\|_{L^2}+\|\nabla\eta\|_{L^2}\|\Delta\eta\|_{L^2}+\|\eta\|_{L^2}\|\Delta\eta\|_{L^2}+\|u\|_{L^2}^{\frac49}\| u\|_{L^5}^{\frac59}\|\nabla\eta\|_{L^2}\|\Delta\eta\|_{L^2}\\
&\leq\frac14\|\Delta u\|_{L^2}^2+\frac14\|\Delta \eta\|_{L^2}^2+C(1+\|u\|_{L^2}^{\frac89}\|u\|_{L^5}^{\frac{10}9})\|\nabla\eta\|_{L^2}^2+C(\|u\|_{L^2}^2+\|\eta\|_{L^2}^2).
\end{align*}
Summing up these estimations, we find that
\begin{align*}
\frac{d}{dt}(\|\nabla u\|_{L^2}^2&+\|\nabla\eta\|_{L^2}^2)+\|\Delta u\|_{L^2}^2+\|\Delta\eta\|_{L^2}^2\\
&\leq C(\|u\|_{L^5}^5+\|u\|_{L^2}^{\frac89}\|u\|_{L^5}^{\frac{10}9})(\|\nabla u\|_{L^2}^2+\|\nabla\eta\|_{L^2}^2)+C(\|\nabla u\|_{L^2}^2+\|\nabla \eta\|_{L^2}^2)+C(\|u\|_{L^2}^2+\|\eta\|_{L^2}^2).
\end{align*}
By applying Gronwall's inequality and using the estimation \eqref{eq30}, we achieve for all $t\leq T^*$
\begin{align*}
\|\nabla u\|_{L_t^\infty(L^2)}^2&+\|\nabla \eta\|_{L_t^\infty(L^2)}^2+\|\Delta u\|_{L_t^2(L^2)}^2+\|\Delta \eta\|_{L_t^2(L^2)}^2\\
&\leq(\|\nabla u_0\|_{L^2}^2+\|\nabla \eta_0\|_{L^2}^2+\beta(T))\exp\left(C\int_0^t\|u\|_{L^5}^5+\|u\|_{L^2}^{\frac89}\|u\|_{L^5}^{\frac{10}9}dt'\right)\\
&\leq(\|\nabla u_0\|_{L^2}^2+\|\nabla \eta_0\|_{L^2}^2+\beta(T))
e^{C_{in}(T)(N^{-1}+N^{-2/9})}\\
&\leq\frac32(\|\nabla u_0\|_{L^2}^2+\|\nabla \eta_0\|_{L^2}^2+\beta(T)),
\end{align*}
provided that $N$ is sufficiently large. $\square$

\section{Acknowledgments}
This work is partially supported by National Natural Science Foundation of China (No.12171111).

\end{document}